\documentclass[1p]{elsarticle}

\usepackage{floatrow}
\usepackage{color} 
    \usepackage{arydshln}
\newtheorem{lemma}{Lemma}

\newcommand{\enma}[1]   {\ensuremath{#1}}

\newcommand{\req}[1]{(\ref{#1.eq})}

\newcommand{\beq}{\begin{equation}}
\newcommand{\eeq}{\end{equation}}
\newcommand{\beqn}{\begin{eqnarray}}
\newcommand{\eeqn}{\end{eqnarray}}
\newcommand{\beqns}{\begin{eqnarray*}}
\newcommand{\eeqns}{\end{eqnarray*}}
\newcommand{\bct}{\begin{center}}
\newcommand{\ect}{\end{center}}
\newcommand{\btmz}{\begin{itemize}}
\newcommand{\etmz}{\end{itemize}}
\newcommand{\benum}{\begin{enumerate}}
\newcommand{\eenum}{\end{enumerate}}

\newcommand{\R}{{\mathbb R}}

\newcommand{\Z}{{\mathbb Z}}

\newcommand{\cO}{\enma{\mathcal O}}

\newcommand{\Imag}[1]{\Im\!\lb #1 \rb}
\newcommand{\Real}[1]{\Re\!\lb #1 \rb}

\newcommand{\obtd}[2]{\matbegin \begin{array}{c:c}
        	#1 & #2    \end{array} \matend }	
\newcommand{\tbod}[2]{\matbegin \begin{array}{c}
        	#1 \\ \hdashline #2    \end{array} \matend }
\newcommand{\tbtd}[4]{\matbegin \begin{array}{c:c}
        	#1 & #2 \\ \hdashline  #3 & #4  \end{array} \matend }

\newcommand{\matbegin}{
    \left[
}
\newcommand{\matend}{
    \right]
}

\newcommand{\bbm}{\begin{bmatrix}} 
\newcommand{\ebm}{\end{bmatrix}} 

\newcommand{\bsm}{\left[ \begin{smallmatrix}} 
\newcommand{\esm}{\end{smallmatrix} \right]} 

\newcommand{\bsbm}{\left[ \begin{smallbmatrix}} 
\newcommand{\esbm}{\end{smallbmatrix} \right]} 

\newcommand{\bbNm}{\begin{bNiceMatrix}} 				
\newcommand{\ebNm}{\end{bNiceMatrix}} 
\newcommand{\bNA}[1]{ \left[ \begin{NiceArray}{#1} } 		
\newcommand{\eNA}{ \end{NiceArray} \right] }

\newcommand{\lb}{\left(}
\newcommand{\rb}{\right)}
\newcommand{\lcb}{\left\{}
\newcommand{\rcb}{\right\}}

\newcommand{\lnb}{\left\|}
\newcommand{\rnb}{\right\|}

\newcommand{\cE}{{\mathcal E}}

\newcommand{\be}{\begin{equation}}
\newcommand{\ee}{\end{equation}}

\newcommand{\cplxs}{ C\kern -.35em \rule{0.03 em}{.7 ex}~   }

\def\complex{\hbox{C\kern -.45em \rule{0.03 em}{1.5 ex}}~}

\newcommand{\wtl}{\tilde{w}}

\newcommand{\rmc}{{\rm c}}

\newcommand{\rmr}{{\rm r}}

\newcommand{\bi}{\begin{itemize}}
\newcommand{\ei}{\end{itemize}}
\newcommand{\ben}{\begin{enumerate}}
\newcommand{\een}{\end{enumerate}}

\newcommand{\cA}{\mathcal{A}}

\newcommand{\cC}{\mathcal{C}}

\newcommand{\bseq}{\begin{subequations}}
\newcommand{\eseq}{\end{subequations}}

\newcommand{\ba}{\begin{array}}
\newcommand{\ea}{\end{array}}

\definecolor{dred}{rgb}{.8,0,0}

\newcommand{\sm}{\text{-}}

\def\clap#1{\hbox to 0pt{\hss#1\hss}}

\newcommand{\btc}{\begin{tabular}{c}}
\newcommand{\btbl}{\begin{tabular}{l}}
\newcommand{\et}{\end{tabular}}

\newcommand{\fs}{\footnotesize}

    \newcommand{\vt}{{\tilde{v}}}

    \newcommand{\zd}{{\dot{z}}}
    \newcommand{\yd}{{\dot{y}}}

	\newcommand{\rom}{\rule{0em}{1em}}

\newcommand{\Ims}[1]{{\sf Im}\!\left( #1 \right) } 
\newcommand{\Nus}[1]{{\sf Nu}\!\left( #1 \right) }

\newcommand{\hsom}{\hspace{1em}} 
\newcommand{\hstm}{\hspace{2em}}

\newcommand{\Span}[1]{ {\sf span} \!\left\{ #1 \right\} }

	\newcommand{\bbms}{\begin{bsmallmatrix}}
	\newcommand{\ebms}{\end{bsmallmatrix}}

	\newenvironment{bmsm}{\renewcommand*{\arraystretch}{.7} \begin{bmatrix}}{\end{bmatrix}}
	\newcommand{\bbmsm}{\begin{bmsm}}
	\newcommand{\ebmsm}{\end{bmsm}}
	
	\newcommand{\rmDx}{{\rm D_x}}
	\newcommand{\rmAr}{{\rm A_n}}
	\newcommand{\rmAk}{{\rm A_k}}
	\newcommand{\rmAz}{{\rm A_0}}
	
	\newcommand{\ML}{M^{\sm \scriptscriptstyle L}}

\usepackage{algorithm}
\usepackage{multirow}

 \usepackage{xcolor} \usepackage{ulem} \usepackage{changebar}

\newfloatcommand{capbtabbox}{table}[][\FBwidth]
\usepackage[utf8]{inputenc}
\usepackage[T1]{fontenc}
\usepackage{textcomp}
\usepackage{amsmath, amssymb,bm}
\usepackage{graphicx}
\usepackage{mathtools}
\usepackage{soul}
\usepackage{amsmath} 

\usepackage{import}
\usepackage{subcaption}
\usepackage{xifthen}
\usepackage{pdfpages}
\usepackage{transparent}

\begin{document}

\begin{frontmatter}

\title{Implicit Boundary Conditions in Partial Differential Equations Discretizations: Identifying Spurious Modes and Model Reduction}

\author[1]{Pascal R. Karam\corref{cor1}}
\ead{p_r_a@ucsb.edu}
\author[1]{Bassam Bamieh}
\ead{bamieh@ucsb.edu}
\affiliation[1]{organization={Department of Mechanical Engineering, University of California at Santa Barbara},
	city={Santa Barbara,CA},
	postcode={93106},
	country={United States}
}
\cortext[cor1]{Corresponding author}
\begin{abstract} 
We revisit the problem of spurious modes that are sometimes encountered  in partial differential equations discretizations. 
It is generally suspected that one of the causes for spurious modes is due to how  boundary conditions are treated, and we use this 
as the starting point of our investigations. 
By regarding boundary conditions as algebraic constraints on a differential equation, we point out that any differential equation 
with homogeneous boundary conditions also admits a typically infinite number of hidden or implicit boundary conditions. 
In most discretization schemes, these additional implicit boundary conditions are violated, and we argue that this is what leads 
to the emergence of spurious modes. 
These observations motivate two definitions of the quality of computed eigenvalues based on violations of derivatives of boundary 
conditions on the one hand, and on the Grassmann distance between subspaces associated with computed eigenspaces on the other. 
Both of these
tests are based on a standardized treatment of boundary conditions and do not require a priori knowledge of eigenvalue locations. 
The effectiveness of these tests is demonstrated on several examples known to have spurious modes.  In addition, these quality 
tests show that in most problems, about half of  the computed spectrum of a differential operator is of low quality. The tests 
also specifically  identify the low accuracy modes, which can then  
be projected out as a type of model reduction scheme. 
\end{abstract}

\end{frontmatter}

\section{Introduction }

Spectral methods and other higher order discretization schemes are attractive due to fast convergence and high accuracy. However in certain problems, these schemes can generate eigenvalues which are non-physical i.e. they do not converge toward an eigenvalue of the original, infinite dimensional system as the size of the discretization is increased. Instead, these so-called "spurious" modes tend to grow (both in quantity and (or) in magnitude) as the approximation is refined. Thus, identification or removal of these "spurious modes" is essential - especially in time marching situations.   
Considerable effort has been dedicated toward avoiding spurious modes. Researchers have proposed many modifications of classic spectral methods to avoid these spurious modes \cite{Bourne2003,Gardner1989a,Huang1994,Lindsay1992,Model1990,Zebib1987,Xu1997,Straughan1996,Fornberg1990}. Additionally, problem specific techniques \cite{Danilov2019,Schroeder1994,Manning2007,Phillips1993}, filtering methods \cite{Bewley1998,Orzag1971}, and careful modification of boundary conditions \cite{canuto1987boundary} have been proposed.  While the aforementioned methods are effective, they are often taxing to implement and do not generalize well. One of the main difficulties in developing a unified theory is that spurious modes arise in two distinct cases: either from the application of boundary conditions or from algebraic constraints (as in the case of the divergence-free constraint in the incompressible limit of the Navier-Stokes equations). To date, we are not aware of a unifying approach to handle spurious modes regardless of discretization method or problem type.  

In this paper we revisit the problem of spurious modes from a new perspective and consider primarily linear problems for clarity. We adopt the view that  boundary conditions of a partial or ordinary differential equation can be thought of as algebraic (in time) constraints, and thus the combined differential equation and boundary conditions form a Differential Algebraic Equation (DAE). This allows us to treat boundary conditions and differential algebraic systems under the same framework.  In section~\ref{BC_Constraints.sec} we show that the interaction of the algebraic constraints and the differential equations always generate additional implicit constraints which the solutions must satisfy. In infinite dimensions, there are typically an infinite number of  implicit (i.e. in addition to the explicitly stated) constraints. 
Most discretization schemes, while enforcing the explicitly given boundary conditions with techniques such as boundary bordering, will fail to enforce those additional implicit constraints. We argue that it is this phenomenon that is responsible for the emergence of spurious modes. 

In Section~\ref{identifying.sec} we show with a specific problematic  example from~\cite{canuto1987boundary} that explicitly enforcing some of the additional 
implicit boundary conditions can eliminate spurious modes. However, it can not be known apriori how many of those additional constraints
should be imposed in any particular discretization. We therefore use the insight gained from this analysis to develop  quality measures 
for any computed eigenvalue/eigenfunction. The first quality measure is termed the ''derivative test'' which quantifies how much the approximate eigenfunctions violate the first implicit constraint.  In numerical experiments, it appears  that this one additional check is sufficient to identify not only classically spurious modes, but also ones that are badly approximated by the discretization (in examples where the actual eigenvalues can be computed analytically). However, 
the threshold of quality measure for discriminating between well and badly approximated eigenvalues remains a user-selected quantity. We develop an alternative quality measure based on the Grassmann distance between eigenspaces (termed the ''angle criterion'') that appears to be more effective in the sense of not requiring  as much user input in determining the threshold level. In addition, this test is applicable to generalized eigenvalue problems.

Lastly, through extensive numerical experiments we find that in any eigenvalue calculation of a discretized differential operator, generally at least half the computed spectrum is of ''low quality''. Since our quality measures quantify the quality of any specific eigenvalue, we can use these measures to project out that part of the spectrum, thus  leading to a model reduction scheme when the differential operator is the infinitesimal generator of a dynamical system. An example presented in Section~\ref{model_reduction.sec} illustrates this kind of model reduction on a wave-like system.

\section{Motivation: Implicit Boundary Conditions and their Approximations} 				\label{BC_Constraints.sec}
We begin with a simple motivating example. Consider the diffusion equation in one spatial dimension with 
	homogeneous Dirichlet boundary conditions
\begin{align}
            	\partial_t ~\psi(x,t) 	 ~&=~  \partial_x^2 ~\psi(x,t), 				& x\in[-1,1], 	~t\geq 0 	,	
																		\label{heat_DIR_ho_PDE.eq}	\\
            					0~&=~ 	\left.	 \psi(x,t)\right|_{x=\pm1} ,	&t\geq 0	.
 																		\label{heat_DIR_ho_BC.eq}
\end{align}
This equation requires an initial (in time) condition $\psi(x,0)={\psi_o}(x)$, and for a very large class of such initial conditions ${\psi_o}$, it has 
infinitely smooth (in both $t$ and $x$) solutions. 

The PDE above  is given with one set of boundary conditions, however, it has an infinite set of ``implicit'' boundary 
conditions. This can be easily seen as follows. Since the  condition $0= 	\left.	 \psi(x,t)\right|_{x=\pm1}$ has to hold for all time, 
the time derivative has to also be zero for all time, and similarly for time derivatives of all orders 
\begin{align} 
	\forall t\geq 0 , ~	\left.	 \psi(x,t)\right|_{x=\pm1} =0
	\hstm &\Rightarrow \hstm 
	\forall t\geq 0 , ~k\geq 0, ~	\partial_t^k \left.	 \psi(x,t)\right|_{x=\pm1} = 0 		\nonumber		\\
		 &\Rightarrow \hstm 
	\forall t\geq 0 , ~k\geq 0, ~	\partial_x^{2k} \left.	 \psi(x,t)\right|_{x=\pm1} =0,			\label{add_BC_heat.eq}
\end{align} 
where the second implication follows from substituting spatial derivatives for time derivative using the PDE~\req{heat_DIR_ho_PDE}.
Thus there is an infinite set of boundary conditions that are ``implicit'', i.e. they arise from a combination of the explicitly stated 
boundary conditions~\req{heat_DIR_ho_BC} together with the PDE~\req{heat_DIR_ho_PDE} itself.

When a numerical scheme is used to discretize the spatial Laplacian in~\req{heat_DIR_ho_PDE} (e.g. finite-differences, Chebyshev 
collocation, Galerkin,  etc.), the approximation is an ODE which can be written in matrix-vector form  
\begin{align}
            		\dot{z}(t) &=~ A ~z(t) , 				& z(t)\in\R^n ,					\label{zd_A.eq}	\\ 
            		0 &=~ C ~z(t) , 						& C \in\R^{q\times n}, ~q<n,		\label{Cz_zero.eq}
\end{align} 
where $n$ is the size of the discretization grid (or the number of basis functions), 
and $q$ is the number of boundary conditions. The   constraint~\req{Cz_zero} 
represents the discretization of the boundary conditions. The system~\req{zd_A}-\req{Cz_zero} is  a
 Differential-Algebraic Equation (DAE), where~\req{Cz_zero} is viewed as an {\em algebraic constraint} since it does not involve any time derivatives, or 
equivalently as a {\em static-in-time} constraint. 
Discretizations of PDEs are usually not written 
 explicitly in the manner shown here, but rather as an unconstrained ODE after boundary conditions are
 incorporated by some technique such as boundary bordering. However, we delay 
the choice of how 
boundary conditions should be handled until we have examined the solution properties of the DAE~\req{zd_A}-\req{Cz_zero}. 

 Consider again a grid-based discretization scheme. Enforcing the constraint~\req{Cz_zero} will enforce the explicitly 
 given boundary conditions~\req{heat_DIR_ho_BC}. However, the approximate solution will likely not satisfy any 
 of the additional implicit boundary conditions~\req{add_BC_heat}. 
 More precisely, the approximations of the even order spatial derivatives
 at the boundaries will not be zero. The diffusion equation is so well behaved that this issue does not cause any perceptible 
 difficulties. 
 
 In this paper we argue that it is the implicit boundary conditions, and the failure of discretizations to approximate them, that cause
 spurious eigenvalues to emerge. To explore this idea, we first have to establish a few properties of DAEs,  we then return to the 
 discussion of PDEs and their discretizations.

\subsection{Finite-Dimensional Setting: Differential Algebraic Equations (DAEs)}

Consider the finite-dimensional 
DAE~\req{zd_A}-\req{Cz_zero}.
%
%
	Similar to the PDE~\req{heat_DIR_ho_PDE}-\req{heat_DIR_ho_BC}, this system also has additional, implicit algebraic 
	constraints. 
	 Indeed, 
	since $0=Cz(t)$ has to hold for all time, then this constraint must be satisfied by all derivatives of $z$, i.e. 
	\be
		\forall t\geq 0, ~0=Cz(t) 
		\hstm \Rightarrow \hstm 
		\forall k\geq 0, ~ 
		0 = Cz^{(k)}(t) .
	  \label{k_der_const.eq} 
	\ee
	From the differential equation~\req{zd_A}, the $k$'th derivative of $z$ can be written in terms of $z$ as 
	\[
		z^{(1)}=A z 
		\hstm \Rightarrow \hstm 
		\left\{ 
		\begin{array}{rcl} 
			z^{(2)} &=&  A z^{(1)} = AA z =  A^2 z 		\\
					& : & 				\\
			z^{(k)} &=& A^k z 
		\end{array} 	\right. ,
	\]
	as can be verified by induction. Combining this expression with the constraints~\req{k_der_const} and rewriting them all 
	together in matrix-vector notation
	\be
		0 =  
		\bbmsm Cz(t) \\ Cz^{(1)}(t) \\ Cz^{(2)}(t) \\ : \ebmsm 
		= 
		\bbmsm Cz(t) \\ CAz(t) \\ CA^2z(t) \\ : \ebmsm 		
		= 
		\bbmsm C \\ CA \\ CA^2 \\ : \ebmsm z(t)  
		\hstm \Leftrightarrow \hstm
		z(t) ~\in~ \Nus{\bbmsm C\\ CA \\ CA^2 \\  : \ebmsm } =:  \Nus{\cO}
	 \label{bigO_def.eq}
	\ee	
	Thus $z(t)$ must always evolve in $\Nus{\cO}$, the null space of the matrix $\cO$ defined above. 
	Note that this null space is contained in the null space of $C$ (i.e. $\Nus{\cO}\subseteq \Nus{C}$), but is 
	generally smaller\footnote{
	In control theory, the subspace $\Nus{\cO}$ is called the {\em unobservable subspace}. The terminology comes from 
		the fact that if $y(t)=Cz(t)$ is regarded as an ``output'' of the system $\dot{z}(t)=Az(t)$, then an initial state in $\Nus{\cO}$ 
		will  
		produce the zero function as an output, and  is thus  indistinguishable from the zero initial state , i.e. it is ``unobservable''
		from the output.}.

		Matrices like $\cO$ in~\req{bigO_def} generated from a pair $(C,A)$ have  special properties. 
		Define the ``truncated'' version $\cO_k$ of $\cO$ as 
		\[
			\cO_k ~:=~  \bbmsm C \\ : \\ CA^{k\sm1}  \ebmsm
		\]
		Two important facts can be established~\cite{hespanha2018linear} about the null spaces of matrices like ${\cO_k}$ 
		as $k$ varies. First, note that as $k$ increases, the null spaces get smaller, but by the 
		 Cayley-Hamilton theorem, the smallest null space is achieved at $k\leq n$, where $n$ is the dimension of $A$, i.e. 
		 \be
            		\begin{aligned}
				\mbox{for any $k,l\geq 0$}, \hstm 
            			k&\geq l 						&
            			\hstm \Rightarrow \hstm 
            			\Nus{\cO_k}&\subseteq \Nus{\cO_l} , 					\\
            			k&\geq n						&
            			\hstm \Rightarrow \hstm 
            			\Nus{\cO_k} &= \Nus{\cO_n}. 
            		  \label{cO_n_def.eq}
            		\end{aligned}
		\ee
		The second fact is that $\Nus{\cO_n}$ is the largest
		 $A$-invariant subspace contained in $\Nus{C}$~\cite{hespanha2018linear}.  
		These facts give a complete characterization of the 
		solvability of the 
		DAE.
	\begin{lemma} 															\label{DAE1.lemma}
		The DAE~\req{zd_A}-\req{Cz_zero} has a solution iff the initial condition is in the  null space $\Nus{\cO_n}$ of the 
		matrix 
		\[
			\def\arraystretch{.7}			
			\cO_n ~:=~  \bbmsm \scriptstyle   C \\ \scriptstyle CA \\ : \\\scriptstyle  CA^{n\sm1} \ebmsm , 
		\]
		where $n$ is the dimension of $z$. 
		When that condition holds, the solution  remains in that  subspace for all $t\geq 0$. 
		$\Nus{\cO_n}$ is the largest $A$-invariant subspace contained in $\Nus{C}$. 
	\end{lemma} 
	\noindent
	Thus although the initial DAE has the algebraic constraint $0=Cz(t)$, there are in fact many other ``hidden'' or implicit 
	constraints that also must hold. The total number of constraints is  $n$ times the number of original explicit constraints 
	$0=Cz(t)$. However, some of these constraints maybe redundant. The total number of non-redundant 
	constraints is the number of linearly independent rows of $\cO_n$, i.e. the rank of $\cO_n$. It may happen that 
	$\cO_n$ is full rank, and in this case the DAE has no feasible solutions for any initial condition other than zero.

	
	Lemma~\ref{DAE1.lemma} states that the solution of the  DAE~\req{zd_A}-\req{Cz_zero} must be confined to a subspace. 
	We can therefore convert the DAE into an unconstrained ODE 
	 by decomposing 
	$\R^n$ into that subspace  and a complement as stated in the following lemma whose proof is in the appendix. 
	\begin{lemma}[\sf Constraint-free Transformation of the DAE] 							\label{DAE2.lemma}
		Consider the following Differential Algebraic Equations (DAE), and the matrix $\cO_n$ constructed from its parameters
            	\be
			\arraycolsep=2pt
            		\begin{array}{rclr}
                        		\dot{z}(t) &=& A ~z(t) , 		&	\hspace{6em}  z(t)\in\R^n ,	\\ 
                        		0 &=& C ~z(t) , 			& C \in\R^{q\times n}, ~q<n,
            		\end{array} 	
			\hspace{6em} 		
			\cO_n ~:=~  \bbmsm \scriptstyle   C \\ \scriptstyle CA \\ : \\\scriptstyle  CA^{n\sm1} \ebmsm . 
            	  \label{lin_DAE_2.eq}
            	\ee
		Let $M$ be a full column rank matrix such  that $\Nus{\cO_n} =\Ims{M}$, and let 
		$\ML$ be any left inverse of it.  Then 
				every solution of the DAE can be written as $z(t)=M y(t)$, 
			where $y$ satisfies  the unconstrained differential equation 
				\be
					\dot{y} (t) ~=~ \big( \ML A M \big) ~y(t)  ~=:~ \rmAr ~y(t) .
				   \label{A_rmno.eq}
			   \ee
	\end{lemma} 
	\noindent
%
%
	

	Lemma~\ref{DAE2.lemma} converts the DAE to an unconstrained differential equation. It is important to note that as a
	dynamical system, the properties of the ``compressed operator'' $\rmAr$ rather than $A$ are important. For example, if one is to characterize 
	stability of the DAE, it is the eigenvalues of $\rmAr$ in~\req{A_rmno} rather than the eigenvalues of $A$ in~\req{lin_DAE_2}
	that characterize stability. 

	\subsection{Implicit/Hidden Boundary Conditions in PDEs} 

	Having seen how algebraic constraints are treated in the finite-dimensional setting, 
%
	we now work through the PDE~\req{heat_DIR_ho_PDE} and~\req{heat_DIR_ho_BC} 
	in analogy with the finite-dimensional case just presented. The  problem is formulated abstractly  as follows 
            	\begin{align}
                        	\partial_t \psi(x,t) ~&=  \partial_x^2 ~\psi(x,t), 	&	
 					 &\Leftrightarrow&
							\tfrac{d}{dt}  \Psi(t) ~&= \cA ~\Psi(t), 					\label{heat_abs_BC.eq}	\\
                       			0~&= \left. 	 \psi(x,t)\right|_{x=\pm1} ,		&
					 &\Leftrightarrow&
							0~&=\cC ~\Psi(t) ,									\label{heat_abs_BC_D.eq}
		\end{align} 	
	where $\Psi(t):= \psi(.,t)$ is the ``spatial profile'' at each time $t$, 
            the operator $\lb\cA \psi\rb(x,t):=\partial_x^2 \psi(x,t)$, and the ``sampling''  operator
            \[
            	\lb \cC  \psi \rb(x,t)  ~:=~ \bbm \psi(1,t) \\ \psi(-1,t) \ebm , 
            \]
          	gives the values of a function  at $\pm1$. 
		Like $\cA$, the sampling operator is unbounded and defined only on a dense subspace. 
		Any PDE with any given boundary conditions can be abstractly represented as above. 

	Just like the finite-dimensional case, the fact that the constraint~\req{heat_abs_BC_D}  has to hold 
	for all time implies  additional implicit boundary conditions. In this case however, there is an infinite number of such 
	 conditions
	\be
                		0 = 
                            	\bbm \cC \\ \cC  \cA \\ : \ebm \Psi (t)
                		=: \cO~ \Psi(t) 
                			\hstm  \Leftrightarrow \hstm 
                		\bbm \psi(\pm1,t) \\  \psi^{(2)}(\pm1,t)  \\ :  \ebm = \bbm 0 \\ 0 \\ : \ebm  , 
    	     \label{obs_op_DN.eq}
            \ee
            The null space  $\Nus{\cO}$ 
             is made up of  functions with  all even order derivatives at the boundaries   zero
             \be
            	\Nus{\cO} ~=~ \lcb \psi:[-1,1]\rightarrow\R; ~~\psi^{(2k)} (\pm1) = 0, ~~k=0,1,\ldots \rcb 	,\label{cX_rmno_d.eq}
            \ee
            This subspace is clearly  $\cA$-invariant since applying $\cA$ increases the order of differentiation by 2. This is 
            precisely the same set of implicit boundary conditions  found earlier in~\req{add_BC_heat}. For a general PDE 
            with homogeneous boundary conditions  written abstractly like~\req{heat_abs_BC}-\req{heat_abs_BC_D}, all 
            the implicit boundary conditions are encoded as the solution being in the subspace 
            \[
            	\Nus{\cO} ~:=~ \Nus{   \bbmsm  \cC \\   \cC \cA \\  : \ebmsm  } . 
            \]

           Now it is important to observe that in any numerical approximation scheme, it is difficult  to ensure (or even approximate) 
           this infinite number of implicit boundary conditions. There is however one exception  to this statement which is useful to 
           keep in mind for contrast with other schemes. If eigenfunctions of the operator $\cA$ with the given boundary condition can be 
           calculated, and they form a basis for the function space over which the problem is defined\footnote{More precisely, 
            	$\lcb v_k \rcb$ should be a Riesz basis.}, then a Galerkin scheme using those eigenfunctions as a basis does indeed
	satisfy the infinite number of implicit boundary conditions as we now explain. 
            
            The  eigenfunctions are calculated from 
             $\cA v=\lambda v$ while imposing boundary conditions. The eigenfunctions however ``automatically'' satisfy 
            all the implicit boundary conditions as well 
            \be
            	\arraycolsep=2pt
		\left. 
            	\begin{array}{rcl} 
            		\cA v &=&  \lambda ~v 	\\
			  0 	&=& \cC ~v
		\end{array}  
		\right\} 
		\hstm\hstm \Rightarrow \hstm \hstm
		\forall k\geq 0, ~~
		\cC \cA^k v ~=~ \cC \lambda^k v ~=~ 0 .
             \label{A_all_der.eq}
            \ee
            Note that this statement is true for any operator and its eigenfunctions since the calculation above did not depend on the 
            specific form of $\cA$ or $\cC$. 
%
%
	Now a Galerkin approximation to a PDE 
            \[
            	\frac{d}{dt} \Psi(t) ~=~ \cA ~\Psi(t) , 
			\hstm \hstm 
			0 ~=~ \cC ~\Psi(t), 
            \]
            using the eigenfunctions up to order $n$ 
          expresses the approximate solution as a finite linear combination 
	 of the basis elements 
	 \be
	 	\Psi(t) ~=~ \sum_{k=1}^n z_k(t) ~v_k , 
	  \label{Psi_Gal.eq}
	 \ee
	 and the differential equations for the time-varying coefficients are simply 
	 \[
	 	\dot{z}_k(t) ~=~ \lambda_k ~z_k(t), 
		\hstm k=1,\ldots, n. 
	 \]
          By~\req{A_all_der}, the eigenfunctions satisfy the infinite number of implicit boundary conditions, and so does the 
          approximate solution~\req{Psi_Gal} since it is a linear combination of such functions. 
            
           The situation outlined above is very special. A Galerkin scheme with functions other than the eigenfunctions of the operator 
           will not in general satisfy all the implicit boundary conditions. Other schemes such as collocation or finite-differences 
           will also not have this very special property.
           Except for the very special situation described above, any finite dimensional approximation of a PDE will not satisfy, 
           and may not even approximate well, any of  
           the infinite number of implicit boundary conditions.

\section{Eliminating and Identifying Spurious Modes}
\label{identifying.sec}
	In the previous section, we argued that the solution to a general DAE satisfies additional "implicit" conditions due to the interaction between
	the algebraic constraint and differential equation. Furthermore, we demonstrated that these implicit constraints are also present in solutions to linear 
	PDEs with homogeneous boundary conditions.  However, when discretizing a PDE operator, one is not guaranteed that the approximate solution satisfies
	any of these implicit constraints. In this section, we demonstrate how this insight enables us to develop tests to determine the "quality" of 
	an approximate eigenvalue/eigenvector. 

	We begin with the following question. 
	Since any finite approximation to a PDE operator is not guaranteed to satisfy the infinite (or even a finite subset) set of constraints  of the original operator,  
	it is natural to ask whether a finite subset of those conditions can be imposed on any such discretization of the original 
	problem. Any linear PDE with $q$ homogeneous boundary conditions and a $n$-dimensional approximation of it can be 
	abstractly expressed as 
	\be
		{\rm PDE}: ~
		\arraycolsep=2pt
		\begin{array}{rcl} 
			\frac{d}{dt} \Psi(t) &=& \cA ~\Psi(t) \\ 
						0 &=& \cC ~\Psi(t) 
		\end{array} 
		\hstm
		\xrightarrow[\mbox{\fs approximation}]{\mbox{\fs finite-dimensional}}
		\hstm
		{\rm DAE}: ~
		\begin{array}{rclcr} 
			\frac{d}{dt} z(t) &=& A ~z(t) ,	& & 	z(t)\in\R^n				\\ 
						0 &=& C ~z(t), 	& & 	C \in\R^{q\times n}, ~q<n.
		\end{array} 		
	  \label{PDE_DAE.eq}
	\ee	
	If we attempt to apply the procedure of Lemma~\ref{DAE2.lemma} to solve the DAE, it turns out that generically 
	$\Nus{\cO_n}=0$, i.e. there is no solution to the  DAE that satisfies all of the (approximations to the) implicit boundary conditions. 

%
		
	Given the above difficulty, we inquire whether explicitly imposing the {\em first few} of the implicit conditions might render 
	an approximate discretization better behaved than just imposing the boundary conditions on their own? As a starting point, we  
	propose the following  procedure for the DAE discretization~\req{PDE_DAE}.
%
%
	\begin{algorithm}
		\caption{}
	\begin{enumerate} 
		\item Choose an integer   $k\geq 1$. 
		\item Form the matrix 
			$			\cO_k ~:=~  \bbm \scriptstyle   C  \\ : \\\scriptstyle  CA^{k\sm1} \ebm . 
			$ 
		\item Find a full-rank  matrix $M$  such that 
				$
                            			\Ims{M}  =   \Nus{\cO_k} 
				$.
		\item 
			Find a left inverse $\ML$ and form
			 the reduced system $\dot{y}(t) = \lb \ML A M \rb y(t) ~=:~ \rmAk~ y(t)$.
%
%
	\end{enumerate} 
	\label{Alg_first}
	\end{algorithm} 
	
	\noindent
	With $k=1$, the procedure above is equivalent to traditional ``boundary bordering'' schemes. For example,  the 
	diffusion equation with homogeneous Dirichlet boundary conditions with $k=1$ 
	gives ${\rm A_1} :=\ML AM$ as the 
	matrix obtained from the Laplacian approximation by removing the first and last rows and columns. 
	
	In Section~\ref{section:implicit_constraints} we apply the above procedure with $k = 1$ to an 
	example problem that produces spurious 
	modes when a Chebyshev collocation scheme is used. 
	Applying the procedure to the same problem with $k>1$ improves the approximation and eliminates the 
	spurious modes, indicating that lack of enforcement of implicit boundary conditions contribute to the 
	emergence of spurious modes. 
	This scheme however has some difficulties which we discuss. Those difficulties lead us to propose that rather than 
	enforce implicit boundary conditions to eliminate spurious modes, it is easier to use those implicit conditions 
	to identify not only spurious modes, but also those that also have relatively low ``approximation quality''. We develop two such 
	quality measures in Section~\ref{eig_tests.sec}.

\subsection{Enforcing Implicit Constraints} 								\label{section:implicit_constraints}					

	The example presented here illustrates how applying Algorithm~\ref{Alg_first}  
	with $\cO_k$,  $k\ll n$ (e.g. with one or two implicit boundary conditions rather than the full set) can significantly 
	improve the properties of a PDE discretization. 
	The following system was studied in~\cite{canuto2007spectral,canuto1987boundary} 
	\be
            \begin{aligned} 
            	\partial_{t}	\bbm  \psi_{1}(t,x) \\ \psi_{2}(t,x) \ebm
			&= - \bbm \tfrac{1}{2}\partial_x & \partial_x \\ \partial_x & \tfrac{1}{2} \partial_x \ebm   
				\bbm \psi_{1}(t,x) \\ \psi_{2}(t,x) \ebm , 
				 \hstm \hstm  && x \in  [-1,1],\; t > 0, 														\\ 
            		\bbm 0 \\ 0 \ebm  &= \bbm \psi_{1}(t,1) \\  \psi_{1}(t,-1) \ebm , && t > 0 .
            \end{aligned}
             \label{canuto_sys.eq}
           \ee
           The eigenvalues of this system are purely imaginary and can be  analytically calculated as 
            \be
	    \bar{\lambda}_k ~=~ i \frac{3 \pi }{8} ~k, 
			\hstm k\in\Z .
             \label{canuto_eigs.eq}
            \ee
            A spatial discretization of this system with a Chebyshev collocation scheme with $n$ grid points for each of the 
            fields $\psi_1$ and $\psi_2$ 
            gives 
            \be 
                    \begin{aligned}
                    	\frac{d }{d t} \begin{bmatrix} z_1(t) \\ z_2(t) \end{bmatrix} 	
        			~&=~   -\begin{bmatrix} \frac{1}{2} & 1 \\ 1 & \frac{1}{2} \end{bmatrix}  
					\begin{bmatrix} \rmDx & 0 \\0 & \rmDx  \end{bmatrix}  
        				 \bbm  z_1(t)  \\ z_2(t)   \ebm 			
        				~=:~ A~ z(t) 									\\ 
                  	\bbm 0 \\ 0 \ebm  ~& =~  \bbm  [0 \cdots 0   1] & 0	\\ [1 0 \cdots 0]  & 0  \ebm 
        				\bbm z_1(t)  \\  z_2(t)  \ebm    
        			~=:~ C ~z(t) 
                    \end{aligned}
            \label{canuto_disc.eq} 
            \ee 
            where $\rmDx\in\R^{n\times n}$ is the first order Chebyshev differentiation matrix, and $z_1$ and $z_2$ are 
            vectors representing the 
	    values of $\psi_1$ and $\psi_2$ respectively at the collocation points. This specific discretization was shown to exhibit spurious (unstable) behavior~\cite[Sec. 4.2.1]{canuto2007spectral}.  The same reference  showed that careful treatment of the PDE at the boundary using characteristic coordinates eliminated these spurious modes.    
	We now show how using Algorithm 1 eliminates spurious behavior without the use of characteristic coordinates.

\begin{figure}[t]
		\centering
		\begin{subfigure}[t]{.32\textwidth}
			\includegraphics[width=\textwidth]{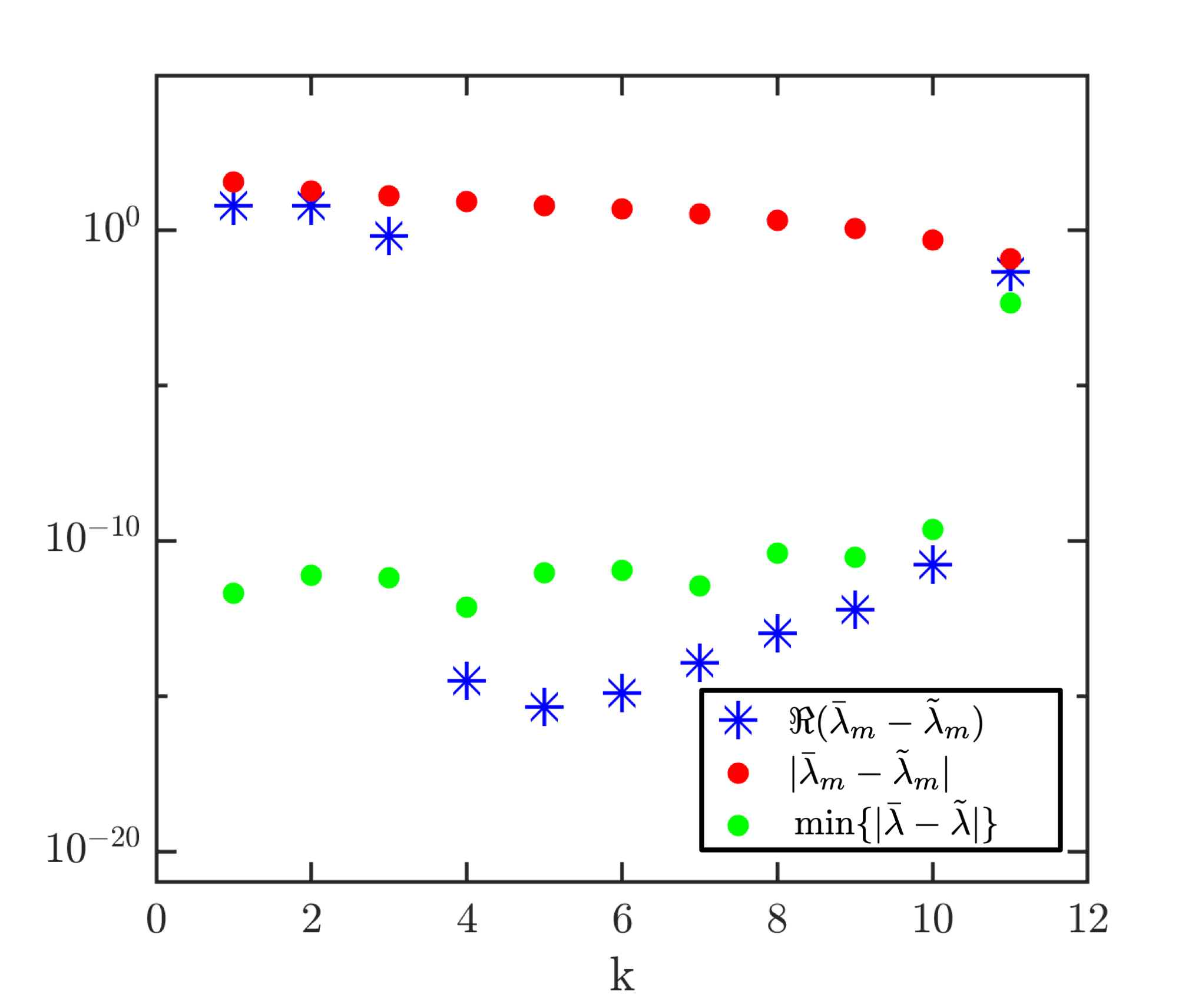}	
			\caption{$N = 16$\label{canuto_k_n16.fig}}
		\end{subfigure}
		\hfill 
		\begin{subfigure}[t]{.32\textwidth}
			\centering 
			\includegraphics[width=\textwidth]{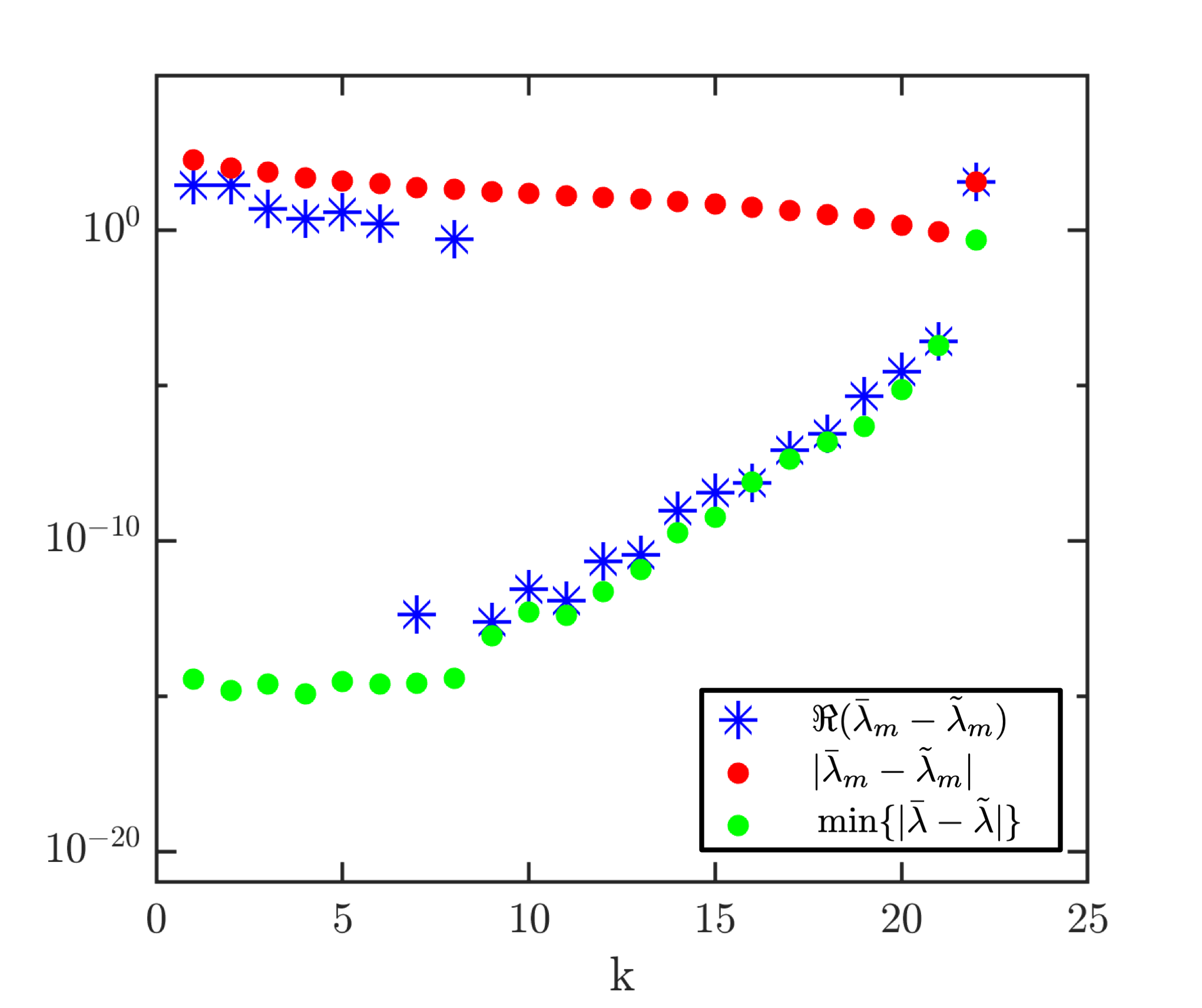}
			\caption{$N = 32$\label{canuto_k_n32.fig}}
		\end{subfigure}
		\hfill
		\begin{subfigure}[t]{.32\textwidth}
			\centering 
			\includegraphics[width=\textwidth]{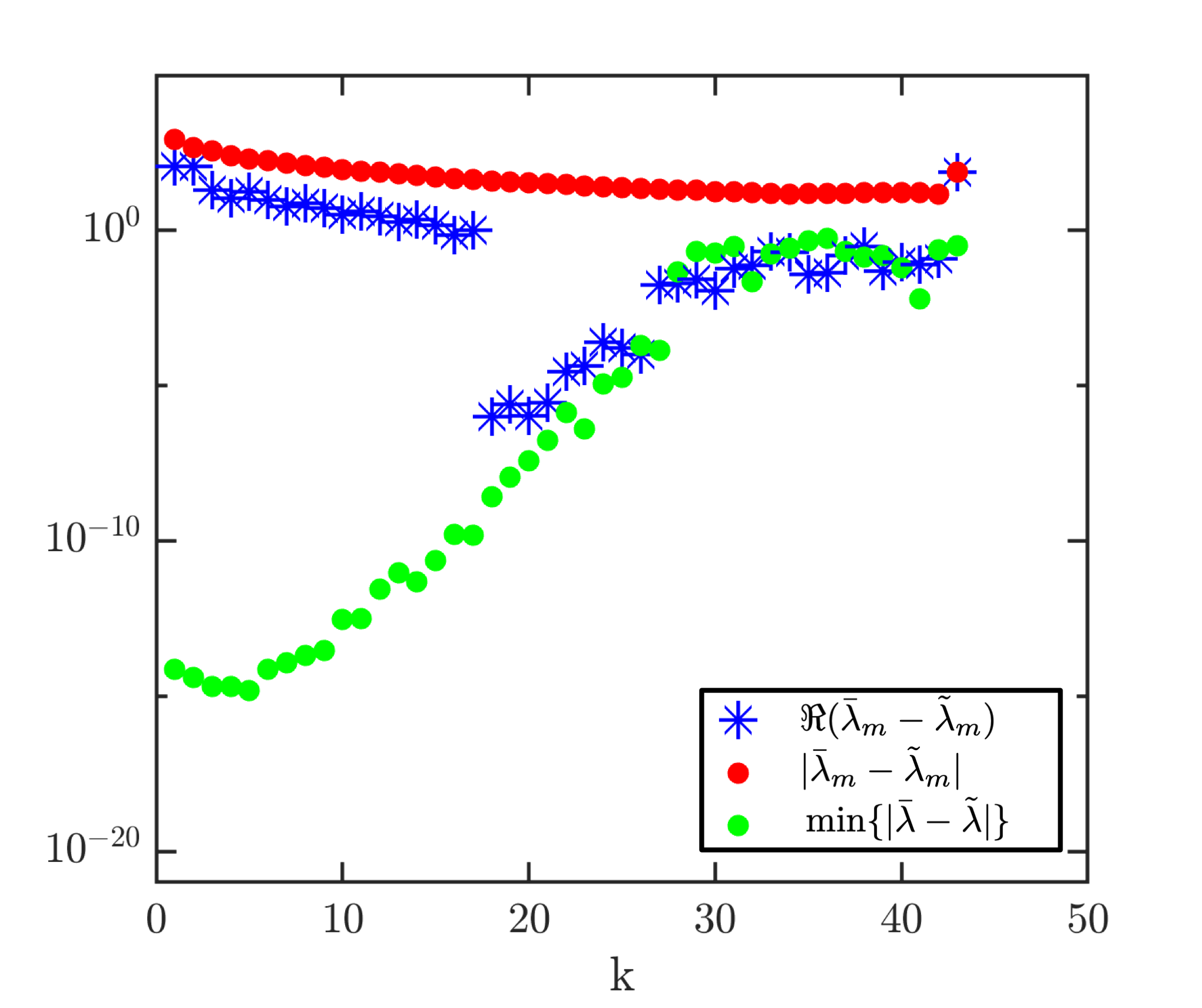}
			
			\caption{$N = 64$\label{canuto_k_n64.fig}}
		\end{subfigure}
		\caption{Explicitly enforcing implicit constraints through $\mathcal{O}_k$ can remove traditionally spurious modes (eigenvalues with erroneously large real part) in spectral collocation schemes. 
			When discretizing  \req{canuto_sys} using a Chebyshev collocation scheme using differentiation matrices of size: $N = 16$ (\subref{canuto_k_n16.fig}), $N = 32$ (\subref{canuto_k_n32.fig}) and $N = 64$ (\subref{canuto_k_n64.fig}) there always exists a sufficiently large $k$ such that the traditionally spurious modes disappear. This is indicated by the sharp drop in real component of the maximum  eigenvalue error. Here $m$ corresponds to the index of the maximum absolute eigenvalue error $m  := \rm{argmax}_i\{ | \tilde{\lambda}_i - \bar{\lambda}_i | \}$    
					However the growth in minimum absolute eigenvalue error for increasing $k$ indicates that applying a large number of implicit constraints contaminates the remaining spectrum. These competing behaviors make it difficult to apply this method in general. 			
			}
		\label{canuto_all.fig}
	\end{figure}

	    The system~\req{canuto_disc} is of the form~\req{lin_DAE_2},  and we can apply Algorithm 1. 
	    For each  $k\geq 1$,  we check the eigenvalues of the resulting 
	reduced matrix $\rmAk = \ML A M$. For this problem, the spurious modes can clearly be identified by their real component as the eigenvalues of the PDE all lie on the imaginary 
	axis by~\req{canuto_eigs}. However, only a small subset of the spectrum of $\rmAk$ has a large real component even though large subset of the spectrum also has a large imaginary error. 
	Therefore, to identify the most problematic eigenvalue, we first find the index of $m$ of the 
	eigenvalue with the largest absolute error as $m:= \underset{i}{\rm{argmax}} | \bar{\lambda}_i - \tilde{\lambda}_i|$ and take the real component of error as 
	\be
	|\Re{(\bar{\lambda}_{m} - \tilde{\lambda}_{m})} |
	 \label{proxy_error.eq}
	\ee
	where $\bar{\lambda}$ denotes the set of eigenvalues \req{canuto_eigs} and $\tilde{\lambda}$ denotes the set of eigenvalues of $\rmAk$. 
	The behavior of \req{proxy_error} as well as the corresponding error's absolute value is shown as a function of $k$ in Figure~\ref{canuto_all.fig} with $N=16$, $N = 32$, and  $N = 64$ collocation point per field. 
	In particular, Figure~\ref{canuto_all.fig} demonstrates that for each discretization size, we can apply a sufficient number of implicit constraints using Algorithm 1 
	such that the spectrum of $\rmAk$ lies purely on the imaginary axis. 
	We find that adding additional implicit constraints into $ \mathcal{O}_{k}$ removes modes with the largest absolute error first which generally correspond to those with large real error. 
	This can be observed by the decreasing trend in absolute eigenvalue error $|\bar{\lambda}_m - \tilde{\lambda}_m|$ in Figure \ref{canuto_all.fig}.
	Therefore, $\rm{dim}(\mathcal{O}_k)$ must be equal or greater than the number of spurious modes with large real real error to guarantee their removal.  
	This is not uniformly true however as see in Figure \ref{canuto_k_n32.fig} where the real component of the error jumps between machine precision and $10^0$ twice in quick succession. 
	
	Despite its potential benefits, we find this approach is not generalizable for two reasons. 
	First, the number of spurious eigenvalues depends on the number collocation points. 
	Thus, without apriori knowledge of the spectrum, it is difficult to identify how many additional constraints are required. 
	This clearly demonstrated in Figure \ref{canuto_k_n32.fig}, \ref{canuto_k_n64.fig} where 9 and 17 additional constraints respectively are required to suppress the spurious real errors. 
	Second, increasing $k$ quickly degrades the quality of the entire spectrum of $\rmAk$ due to the ill-conditioned derivative matrices required to generate $\mathcal{O}_k$.
	This can be observed in Figure \ref{canuto_all.fig} as $\min \{ | \bar{\lambda} - \tilde{\lambda} \}$ grows consistently with $k$.
	Interestingly, in both Figure \ref{canuto_k_n32.fig} and \ref{canuto_k_n64.fig}, the minimum spectral error grows steadily with increasing $k$ after $k \approx 10$.
	This phenomenon is consistent in experiments up to $N = 512$ (not depicted here) which indicates that the minimum spectral error depends only on  $k$ - irrespective of $N$.
	The generally unknown number of spurious modes combined with the degrading performance for even moderate $k$ makes this approach difficult to apply in general. 
	Instead of applying the constraints directly, we now show that checking to see if approximate eigenfunctions satisfy implicit constraints provides a good indicator for mode quality. 
	

	\subsection{Tests for Identifying Spurious Modes} 				\label{eig_tests.sec}

	Enforcing  additional implicit constraints, while effective in some cases, has 
	 the difficulties described  in the previous section. However, we can use the basic fact that violations of implicit constraints 
	 lead to lower quality approximations of a PDEs discretization to develop ``quality certificates'' for numerically 
	 computed eigenvalues. The quality certificates do not require knowledge of the actual eigenvalues.  
	 We develop two different quality measures which both measure the violation of the first implicit constraint: the first calculates the how well eigenfunctions satisfy implicit constraint directly, the second utilizes the notion of principle angles between subspaces to quantify eigenfunction quality.

%
%

\subsubsection{Time Derivatives of Eigenfunctions at Boundaries as a Violation Measure}			\label{der.sec}

The basic idea behind the boundary derivative criterion is to enforce only the given boundary conditions and none of the implicit 
ones, then use the violation of the first implicit boundary condition as an inverse  measure of quality. One can of course 
employ further implicit boundary conditions similarly, but from 
 our numerical experiments, it appears that using only 
the first implicit constraint appears to be sufficient to capture this quality measure. 

Recall again the relation between the original PDE, the  DAE arising from a spatial discretization,  and the unconstrained system after incorporating the given boundary conditions (but not the implicit ones) 
\[
	\arraycolsep=2.5pt
	\begin{array}{rcl} 
		\dot{\Psi}(t) &=& \cA ~\Psi(t) \\ 
		0 &=& \cC ~\Psi(t) 
	\end{array} 
	\hstm 
		\xrightarrow[\mbox{\fs discretization}]{\mbox{\fs spatial}}
	 \hstm 
	\begin{array}{rcl} 
		\zd(t) &=& A ~z(t) \\ 
		0 &=& C ~z(t) 
	\end{array} 
	\hstm 
		\xrightarrow{z(t):=My(t)}
	 \hstm 
	 \begin{array}{rll}
		\yd(t) &=& \lb \ML A M \rb~ y(t) 		\\
		       &=:& \\rm{A_1} ~y(t)
	\end{array}
\]
where the  matrix  $M$ is chosen so that $\Ims{M}=\Nus{C}$. 
The solution $z$ of the DAE is then reconstructed from the unconstrained solution $y$ by $z(t)=My(t)$.  Two points 
deserve pointing out:
\begin{enumerate} 
	\item By construction, the solution $z(t)=My(t)$ satisfies the boundary conditions $0=Cz(t)$ since 
		\[
			C~z(t) ~=~ C M ~y(t) ~=~ 0 		\tag{$\Ims{M}=\Nus{C} ~\Rightarrow~ CM=0$}. 
		\]
	\item The reconstructed solution $z(t)=My(t)$ may not satisfy the   differential equation for $z$ precisely 
		since 
		\be
            				\dot{z}(t) ~=~  M~\yd(t) ~=~ M\rmAz ~y(t) ~=~ M\ML AM~ y(t)  ~=~ M\ML A ~z(t)	
		\label{z_y_ineq.eq}
		\ee
		Note that $M\ML=I$ only on $\Ims{M}$, and for equality above to hold $AMy(t)$ must be in $\Ims{M}$. 
		However, since $\Ims{M}$ is not guaranteed to be $A$-invariant, this may not be the case. Recall that 
		$\Ims{M}\subseteq\Nus{C}$ is guaranteed to be $A$-invariant only if 
		$\Ims{M}=\Nus{\cO_n}$, which is not the case in this scheme. 
\end{enumerate} 
%
%
The observations in~\req{z_y_ineq} motivate several quality measures of the solution $z(t)=My(t)$ obtained from the 
unconstrained problem. For example, over a time interval $[0,T]$, the average (mean square)
 violation of the differential equation 
is 
\[
	\frac{1}{T} \int_0^T \lnb \zd(t) - M \yd(t) \rnb^2  dt 
	~=~ 
	\frac{1}{T} \int_0^T \lnb \lb AM - M\ML A M\rom  \rb  y(t) \rnb^2  dt 
	~=~ 
	\frac{1}{T} \int_0^T \lnb  (I - M\ML) A M  ~  y(t) \rnb^2  dt .
\]
Another quality measure is the average violation of the first implicit boundary condition, i.e. the boundary condition on 
$\zd$ 
\[
	\frac{1}{T} \int_0^T \lnb C~\zd(t)  \rnb^2  dt 
	~=~ 
	\frac{1}{T} \int_0^T \lnb CA~z(t)  \rnb^2  dt 
	~=~ 
	\frac{1}{T} \int_0^T \lnb CAM ~y(t)  \rnb^2  dt .
\]

While the above criteria can measure the quality of a solution over time, their evaluations require actually solving the 
equation for $y$, which is not possible over significant time intervals in the presence of unstable spurious modes. 
In problems where we are mainly interested in computing eigenvalues (even unstable ones), it is more useful 
to use the above criteria to test the quality of any computed eigenvector of $\rm{A_1}$. Let $(\tilde{\lambda},\vt)$ be an 
eigenvalue/vector pair computed for $\rm{A_{1}}$. If we use $z(0)=M\vt$ as an initial condition for the DAE, we expect the 
initial derivative $\zd(0)$ to satisfy the boundary condition $C\zd(0)=0$. However, the quantity 
\[
	C~\zd(0) ~=~ CA~z(0) ~=~ CAM ~\vt 
\]
may not be zero. We can adopt the norm of this quantity as a measure of the violation of the first implicit boundary 
condition by the computed (now thought of as ``proposed'') mode $M\vt$. We state this formally in Algorithm \ref{Alg_second}. 
%
%
%
%
%
%
\begin{algorithm}												
	\caption{Implicit Constraint (Boundary Derivative) Test}		
	\begin{enumerate}
		\item Form the unconstrained, compressed system 
			$\dot{y}(t) = \lb \ML A M\rb  y(t) =: \rm{A_1} \, y(t) $, where $M$ is full column rank
			 with $\Ims{M}=\Nus{C}$
					(this corresponds to the choice $\mathcal{O}_1 = C$ in Algorithm 1).
				\item  Calculate the eigenvectors $\lcb  \tilde{v}_1, \tilde{v}_2, \ldots, \tilde{v}_m \rcb$  of the compressed system $\rm{A_1}$. 
					\item Normalize each eigenvector to have 2-norm 1: $\left\lVert v_k \right\rVert = 1, \; \forall k \in 1,\ldots,m$
		\item	 For each eigenvector $\vt_l$, test whether it satisfies the first implicit constraint by computing   
			 $s_l :=  C A M \tilde{v}_l$.
		\item		$\left\lVert s_l \right\rVert$ 
				quantifies how much $\vt_l$ {\em violates}  the first implicit constraint. We define  $1/\left\lVert s_l \right\rVert$
				as a measure of the quality of that numerically computed eigenvalue $\tilde{\lambda}_l$. 
	\end{enumerate}
\label{Alg_second}
\end{algorithm}

\begin{figure}[t]
		\centering
		\begin{subfigure}[t]{.32\textwidth}
			\includegraphics[width=\textwidth]{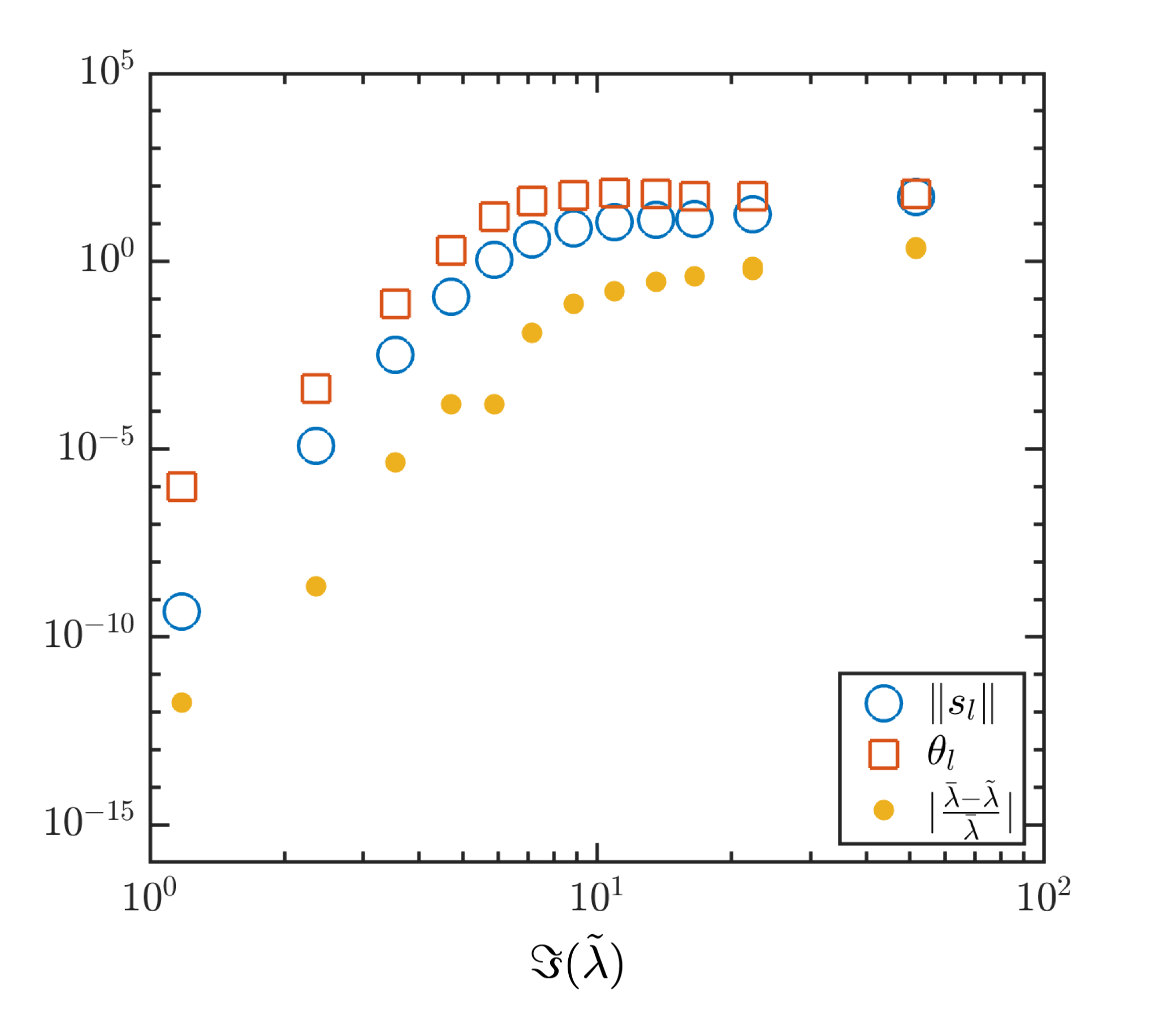}	
			\caption{$N = 16$\label{Canuto_1D_rel_error16.fig}}
		\end{subfigure}
		\hfill 
		\begin{subfigure}[t]{.32\textwidth}
			\centering 
			\includegraphics[width=\textwidth]{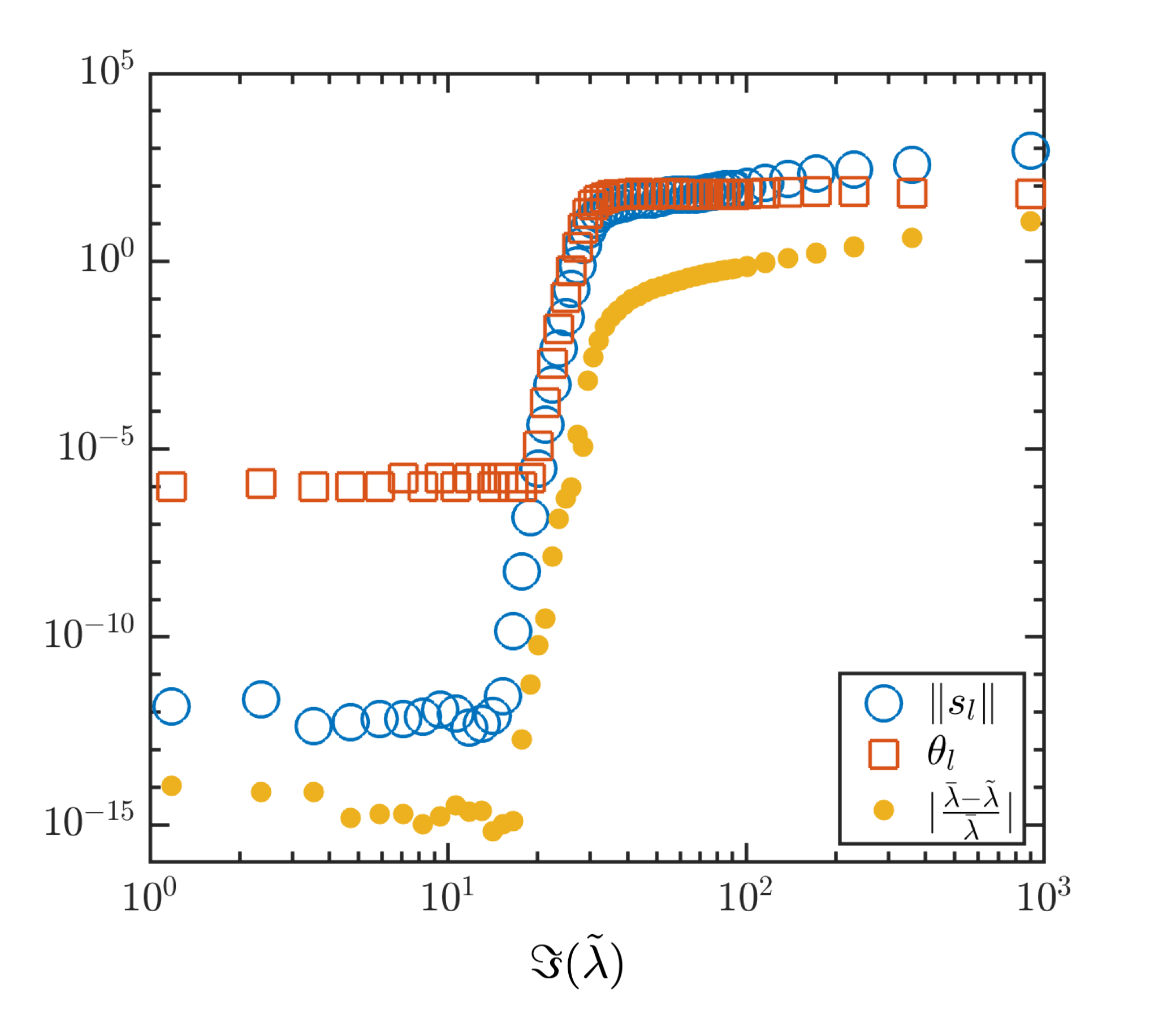}
			\caption{$N = 64$\label{Canuto_1D_rel_error64.fig}}
		\end{subfigure}
		\hfill
		\begin{subfigure}[t]{.32\textwidth}
			\centering 
			\includegraphics[width=\textwidth]{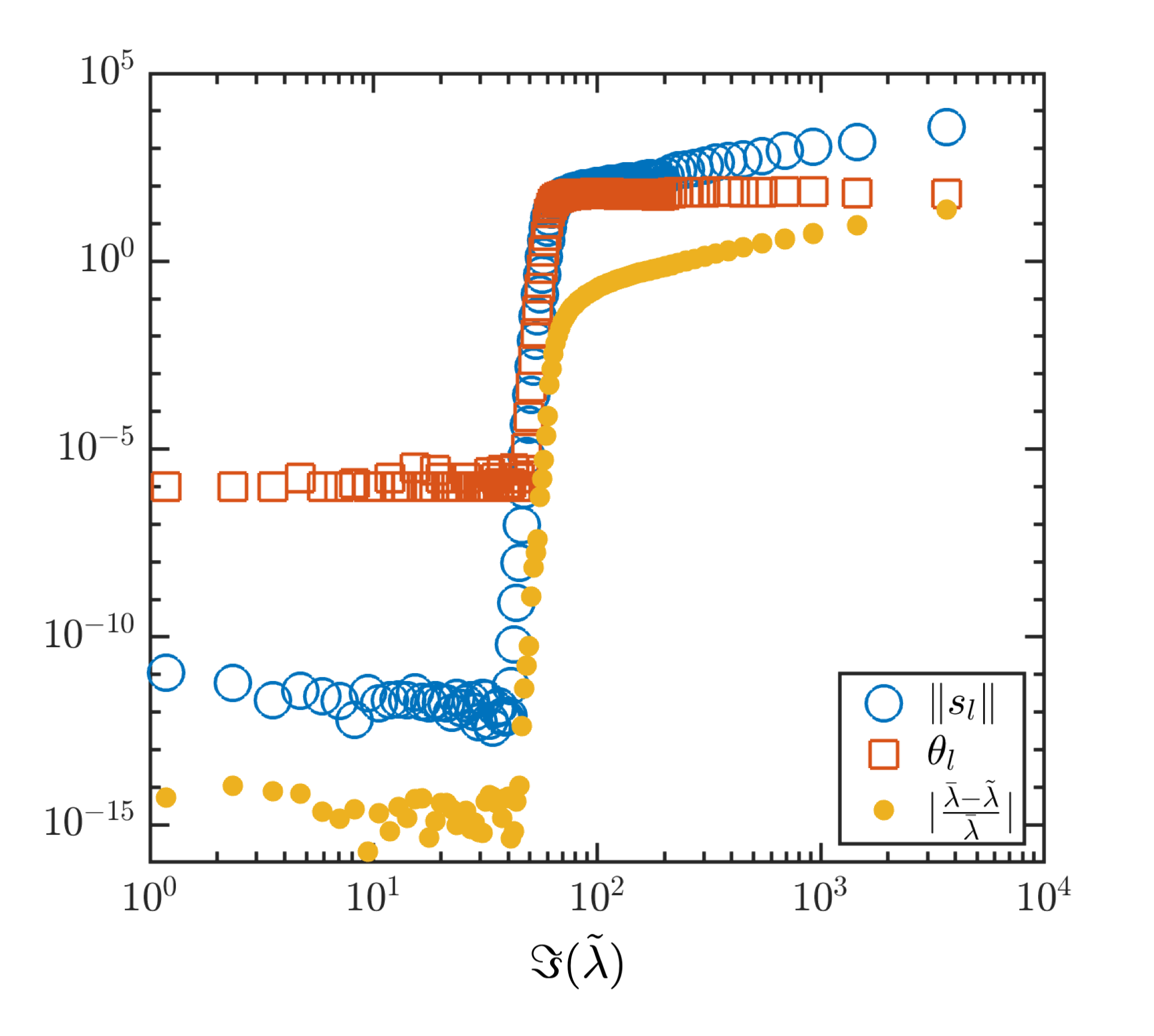}
			
			\caption{$N = 128$\label{Canuto_1D_rel_error128.fig}}
		\end{subfigure}
	\caption{
			Analysis of the ``quality'' of the spectrum of system~\req{canuto_disc} 
			 discretized using a Chebyshev collocation scheme with $N = 16$ (\subref{Canuto_1D_rel_error16.fig}), $N = 64$ (\subref{Canuto_1D_rel_error64.fig}), and $N = 128$ (\subref{Canuto_1D_rel_error128.fig})  nodes per field. For each discretization size, the three error criteria are plotted versus the imaginary component of the respective eigenvalue. 
			 The relative error $(\bar{\lambda}-\tilde{\lambda})/\bar{\lambda}$, the first implicit constraint violation measure $\|s_l\|$  (derivative test), and the Grassmann distance $\theta_l$ are shown for all three cases. 
			 Boundary conditions are implemented using the automated procedure (Step 1) of Algorithm~\ref{Alg_second}. 
			 In all three cases, both error criteria the relative eigenvalue error well. 
			 Note that because the Grassmann distance is root sum square of two angles, a Grassmann distance of  $1e \sm 6$ equates to an error of order machine precision. 
			} 
			\label{numvstheo.fig}
\end{figure}

Our numerical experiments on the example~\req{canuto_sys}-\req{canuto_disc} show that $\left\lVert s_l \right\rVert$ tracks the relative eigenvalue error $\lb{\bar{\lambda}_l - \tilde{\lambda}_l}\rb/{\bar{\lambda}_l}$ closely. 
To demonstrate we return to \req{canuto_disc} approximated using a Chebyshev collocation scheme with $N = 16$,  $N = 64$, $N = 128$ collocation points per field. 
We find that boundary derivative criterion $\left\lVert s_{l} \right\rVert$ tracks the relative eigenvalue  error closely $|\frac{\bar{\lambda} - \tilde{\lambda}}{\bar{\lambda}}|$ (Figure \ref{numvstheo.fig}). 
Furthermore, these experiments demonstrate that even though a minority of the spectrum of $\rm{A_1}$ is classically spurious (has large real component), a significant part of the spectrum has errors in the imaginary component. 
While these errors are not as detrimental, as they does not alter stability, they do not accurately reflect the behavior of the original PDE.  

Though results are promising for simple cases such as \req{canuto_disc}, Algorithm 2 comes with two main drawbacks.
First, it cannot be expanded to handle partial differential algebraic equations (PDAEs). Finding the spectral decomposition of PDAEs amounts to solving a generalized eigenvalue problem.  In problems of this form, one cannot obtain an expression for the derivative in time at the boundary explicitly in terms  of the eigenfunctions.   
In addition, evaluating $\left\lVert s_l \right\rVert$ for systems with boundary conditions containing higher order  derivatives can be questionable  due the ill-conditioned higher derivative matrices.

\subsubsection{Angle Criterion for Eigenvalue Problems}					  \label{angle.sec}

For a differential operator with boundary conditions, the eigenvalue problems for the original operator, its spatial 
discretization, and the ``compressed'' (unconstrained) system are 
\[
	\arraycolsep=2.5pt
	\begin{array}{rcl} 
		\cA ~\psi &=& \lambda~ \psi   \\ 
		0 &=& \cC ~\psi 
	\end{array} 
	\hstm 
		\xrightarrow[\mbox{\fs discretization}]{\mbox{\fs spatial}}
	 \hstm 
	\begin{array}{rcl} 
		A ~w  &=&  \lambda ~w		\\
		0 &=& C ~w 
	\end{array} 
	\hstm 
		\xrightarrow{w:=Mv}
	 \hstm 
	 \begin{array}{rll}
		 \rm{A_1} ~v  &=& \lambda ~ v 
	\end{array}, 
\]
where again the matrix $M$ is chosen so that $\Nus{C}=\Ims{M}$. Once an eigenvector $\vt$  of $\rm{A_1}$ is found (say 
$A\vt=\tilde{\lambda} \vt$), then 
the ``proposed'' vector $\wtl:=M\vt$ should satisfy the equations $A\wtl=\tilde{\lambda} \wtl$ and $C\wtl =0$. 
The constraint $C\wtl=CM\vt = 0$ is satisfied by construction. However, the equation $A\wtl=\tilde{\lambda} \wtl$
is not guaranteed to hold. This motivates another error criterion quantifying how badly this equation is violated. 

Observe that $A$ is the discretization of the differential operator $\cA$ without regard to boundary conditions. An 
eigenfunction $\psi$  of $\cA$ that satisfies the boundary conditions must satisfy both equations $\cA\psi=\lambda \psi$ 
and $\cC \psi =0$. Since $\wtl=M\vt$ is the numerically computed approximation to an eigenfunction $\psi$, it should 
also satisfy the equation $A\wtl = \tilde{\lambda} \wtl$, which is the discretization of the eigenfunction equation before 
incorporating boundary conditions. 
Since eigenvectors are only defined up to multiplication by complex scalars, we need to quantify the violation of the 
 equation  $A\wtl = \tilde{\lambda} \wtl$ where $\wtl$ is replaced by an entire subspace rather than a single vector. 
 The Grassmann distance computed from the principal angles between subspaces provides such an error measure. 

Complex
eigenvectors (with geometric multiplicity one)
 are  defined only up to a one-dimensional complex subspace, or equivalently a two-dimensional real subspace. 
 Thus the 
criterion for $\tilde{w}=M\vt$ to be an eigenvector for $A$ is 
			  \[
				\Span{\Real{AM\vt},\Imag{AM\vt}} 
				~\subseteq~ 
			  	\Span{\Real{M\vt}, \Imag{M\vt}}.
			  \]
If the eigenvalue has geometric multiplicity $\alpha$, and $\lcb \vt_1,\ldots,\vt_\alpha \rcb$ are a basis for its invariant subspace, then 
the criterion would become		
			  \[
				\Span{\Real{AM\vt_1},\Imag{AM\vt_1},..,\Real{AM\vt_\alpha},\Imag{AM\vt_\alpha}} 
				~\subseteq~
			  	\Span{\Real{M\vt_1}, \Imag{M\vt_1},..,\Real{M\vt_\alpha}, \Imag{M\vt_\alpha}}.
			  \]
Thus testing whether $\tilde{w}=M\vt$ is indeed an eigenvector of $A$ amounts to testing whether two subspaces coincide. 			  
The most common measure of distances between subspaces is the 
			  {\em Grassmann distance} computed from 
			  the principal angles  between  two subspaces. 
			  We will use this distance as a quality measure of the computed eigenvalues and show with 
			  a collection of examples that this criterion is highly effective in identifying spurious modes. 
			  An additional feature of this criterion is that it is also applicable  to
			    generalized
			  eigenvalue problems.

For simplicity we only consider examples with (complex) geometric multiplicity of one. This means calculating the Grassmann 
distance between two 2-dimensional real subspaces which we now briefly summarize. 
 Let $u_1,u_2 \in \mathbb{C}^{N} $ be vectors, and consider the two 
real subspaces   
\[
	{\cal U}_1 := \Span{ \Real{u_1}, \Imag{u_1}} ~\subset ~ \R^{N}, 
	\hstm\hstm 
	{\cal U}_2 := \Span{ \Real{u_2}, \Imag{u_2}} ~\subset ~ \R^{N} .
\]
Let $U_1,U_2$ be matrices whose columns are orthonormal bases for ${\cal U}_1,{\cal U}_2$ respectively. 
The so-called principle angles $\theta_1,\theta_2$ between the two subspaces are calculated  from the singular values of 
the $2\times 2$ matrix $U_1^*U_2$ by 
\begin{equation}
	\mbox{singular values}\lb \rom U_1^{*} U_2 \rb 
		 = \lcb  \cos\theta_1 , \cos \theta_2 \rom \rcb. 
\end{equation}
The Grassmann distance is then $\theta = \sqrt{\theta_1^2 + \theta_2^2} $, a quantity which in this case is between $0$ and 
$90\sqrt{2} $. This leads to Algorithm~\ref{Algorithm_third} for quantifying  a computed eigenvalue quality. 

\begin{algorithm}
	\caption{Grassmann Distance Test (Angle Criterion)}
\begin{enumerate}
	\item Form the compressed system $\rm{A_1} := \ML A M$, where $M$ is full column rank with $\Ims{M}=\Nus{C}$.
	\item  Calculate the eigenvectors $\lcb  \tilde{v}_1, \tilde{v}_2, \ldots, \tilde{v}_m \rcb$  of the compressed system $\rm{A_1}$. 
		\item Normalize each eigenvector to have 2-norm 1: $\left\lVert v_k \right\rVert = 1, \; \forall k \in 1,\ldots,m$
	\item For each eigenvector $\tilde{v}_l$, calculate the Grassmann distance $\theta_l$ between $ M \tilde{v}_l$ and $A M \tilde{v}_l$.
	\item $\theta_l$ quantifies how ``incongruent'' the two subspaces are.  We define  $1/\theta_l$
				as a measure of the quality of that numerically computed eigenvalue $\tilde{\lambda}_l$. 
\end{enumerate}
\label{Algorithm_third}
\end{algorithm}

Our experiments indicate that $\theta_l$ generally provides slightly better results for eigenvalue/eigenvector quality than 
the derivative test $\left\lVert s_l \right\rVert$ of the previous section. Specifically, we find that $\theta_l$ tracks the relative eigenvalue error consistently regardless of the size of the discretization i.e. is relatively robust against round-off error (See example in Section \ref{section:generalized} and Table~\ref{Orr_Som.table}).  Conversely, for the same eigenvalue, $\left\lVert s_l \right\rVert$ tends to grow with increasing discretization size even if the relative eigenvalue error remains constant . This makes it difficult to implement consistent thresholds to identify well captured eigenfunctions. However, $\theta_l$ does not come without drawbacks. In practice, care must be taken with eigenfunctions corresponding to zero eigenvalues.  When calculating $A M \tilde{v}$ corresponding to an eigenvalue at zero, even if the eigenfunction is well captured, $A M \tilde{v}  = O(10^{-16})$. The remaining noise will result in an erroneously large Grassmann distance. For eigenfunctions corresponding to eigenvalues at zero, we recommend simply checking the size of $A M  \tilde{v}$. If $\left\lVert A M \tilde{v} \right\rVert_{2} = O(10^{-16})$, then $M \tilde{v}$ is well captured.

%

\subsubsection{Generalized Eigenvalue Problems}
\label{section:generalized}

The Grassmann Distance test can easily be expanded to handle generalized eigenvalue problems which arise when working with partial differential  equations of the form~\cite{Manning2007}
\be
	\begin{aligned}
		\partial_t ~\cE \psi ~&=~ \cA \psi , 			\\ 
					0 ~&=~ \cC \psi 
	\end{aligned} 
  \label{descriptor.eq}
\ee
where $\cE,\cA$ are  spatial PDE operators, and the operator $\cC$ encodes homogeneous boundary conditions as well as other algebraic constraints. 
For example, the the Orr-Sommerfeld equation  for plane Poiseuille flow in two dimensions ($z$ as the wall-normal direction and $x$ as the 
streamwise direction)
\be
	\begin{aligned}
	\partial_t ~ \alpha R \lb \partial_{z} ^{2} - \alpha ^{2} \rb ~~\psi_\alpha(z,t)  
	~&=~
	\lb \partial_{z}^{4} + \lb -2 \alpha ^{2} - \alpha R i \bar{u}(z) \rb ~ \partial_{z}^{2} 
			+ \lb  \alpha ^{4} + \alpha ^{3} R i \bar{u}(z) - 2R \alpha i \rb \rom \rb~~ \psi_\alpha( z, t) 					\\
						  \bbm 0 \\ 0 \ebm  ~&= ~ \bbm  \psi_\alpha(\pm 1,t) \\  \partial_z \psi_\alpha (\pm 1,t) \ebm , 
	\end{aligned}
  \label{OS_PDE.eq}
\ee
where $\alpha$ is the streamwise wavenumber, $R$ is the Reynolds number and $\bar{u}(z) = 1-z^{2}$ is the laminar Poiseuille flow profile.
Note that $\psi$ 
here is the stream function after Fourier transformation in the streamwise direction (which replaces the streamwise coordinate $x$ with 
the wavenumber $\alpha$), and~\req{OS_PDE} represents a parameterized (by $\alpha$) family of PDEs in one spatial dimension.

The generalized eigenvalue problem for systems of the form~\req{descriptor} correspond to finding eigenfunctions $\psi$ that satisfy 
\[
	\lb \lambda \cE - \cA \rb \psi ~=~ 0
	\hstm \Leftrightarrow \hstm 
	 \lambda ~ \cE \psi = \cA \psi  .
\]
Let $A$ and $E$ be discretizations of the spatial PDE operators $\cA$ and $\cE$ respectively (before incorporation of any boundary 
conditions). Let  $\wtl$ be any candidate discretization of an eigenfunction, then the angle criterion we adopt is the Grassmann distance 
between the two subspaces 
\[
	\Span{\Real{E\wtl},\Imag{E\wtl}} ,
		\hstm 
	\Span{\Real{A\wtl}, \Imag{A\wtl}}.
\]
The vector $\wtl$ can come from any type of calculation. For example, $\wtl=M\vt$ where $M$ has orthonormal columns and  $\Ims{M}=\Nus{C}$, and the vector $\vt$ 
calculated from the the generalized (compressed) eigenvalue problem 
\[
	\tilde{\lambda} ~\ML EM ~\vt ~=~ \ML AM ~\vt.
\]

%

\begin{figure}[t]
	\begin{floatrow}
		\ffigbox{\includegraphics[width=0.5\textwidth]{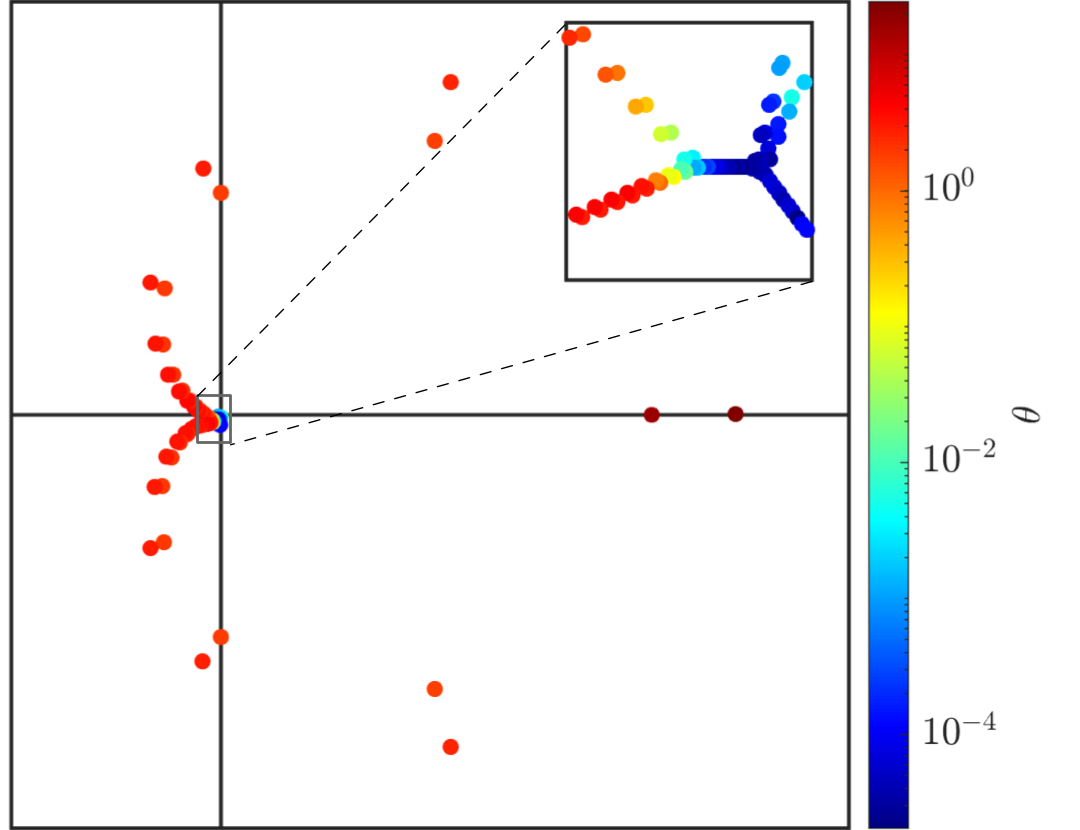}}{
	\caption{Spectrum of Orr-Sommerfeld operator for Poiseuille flow with R = 10, 000, $\alpha$ = 1. The system is discretized using a Chebyshev collocation scheme with a grid size of N = 130. The approximate system contains a large number of spurious modes which are clearly identified by the Grassmann distance criterion $\theta$.  Inset shows the classic Orr-Sommerfeld spectrum which is captured well by the discretized problem.}
	\label{fig:Orr-Som_spectrum}
}
\capbtabbox{
	\begin{tabular}{||r|r|l||}
		\hline
		\multicolumn{1}{||c}{$\theta$ } & \multicolumn{1}{|c|}{$|\tilde{\lambda} - \bar{\lambda}|$} & \multicolumn{1}{c||}{N} \\ 
\hline 
\hline
4.498e-05 & 8.956e-07 & 80 \\ 
3.617e-04 & 3.611e-06 & 90 \\ 
1.001e-03 & 8.801e-06 & 100 \\ 
2.498e-03 & 2.378e-05 & 110 \\ 
4.931e-03 & 1.185e-04 & 70 \\ 
6.622e-03 & 5.539e-05 & 120 \\ 
8.785e-03 & 7.941e-05 & 140 \\ 
1.944e-02 & 1.331e-04 & 130 \\ 
1.173e-01 & 1.081e-03 & 150 \\ 
4.013e-01 & 3.720e-03 & 60 \\ 
1.541e+00 & 4.047e-02 & 50 \\ 
\hline 
	\end{tabular}}{\caption{Approximate eigenvalue $\tilde{\lambda}$ corresponding to the Tollmien-Schlichting (T-S) wave of the Orr-Sommerfeld operator for Poiseuille flow discretized using Chebyshev collocation with R = 10,000, $\alpha = 1$, grid sizes N = $50,\ldots, 150$. This is the same configuration as in \cite{Orzag1971} who reported the eigenvalue corresponding to the T-S wave to be $0.00373967 - 0.23752649i$. Our reported eigenvalue with smallest Grassmann distance $\theta$ corresponds best with that in \cite{Orzag1971}. }
\label{Orr_Som.table}}
\end{floatrow}
\end{figure}

A well-studied~\cite{Orzag1971} instance  of this equation corresponds to $\alpha = 1$ and  $R =10,000$, where the system is known 
to have an instability. 
In addition, the eigenvalues of this system are notoriously sensitive to small variations in problem parameters, which is 
one reason it has become a benchmark test case for numerical algorithms. 
The generalized eigenvalue problem in this case becomes 
\begin{align}
	\lambda~ \alpha R \lb \partial_{z} ^{2} - \alpha ^{2} \rb \psi(z) 
	~&=
	\lb \partial_{z}^{4} + \lb -2 \alpha ^{2} - \alpha R i \bar{u}(z) \rb ~ \partial_{z}^{2} 
			+ \lb  \alpha ^{4} + \alpha ^{3} R i \bar{u}(z) - 2R \alpha i \rb \rom \rb \psi(z) ,	
						\hsom z\in[-1,1], 											\nonumber 	 \\ 
						0 ~&=~
						     \psi(\pm 1) = \partial_z \psi (\pm 1) .
  \label{orr-som.eq}
 \end{align}
Generally, when this system is approximated using a Chebyshev collocation method, a large portion of the spectrum is spurious as 
shown in Figure~\ref{fig:Orr-Som_spectrum} where a  a collocation grid of size $ N = 130$ is used. The Grassmann distance criterion clearly identifies the portion of the spectrum which is accurately captured versus the spurious parts. 

Other tests of this example are reported in Table~\ref{Orr_Som.table}. 
With the same values of $\alpha$ and $R$ used above,  Orzag~\cite{Orzag1971} reported the eigenvalue corresponding to the Tollmien-Schlichting (T-S) mode to be $0.00373967 - 0.23752649i$ (real and imaginary components switched to match the configuration used in this paper). 
The table shows the eigenvalue closest to the T-S mode for varying grid sizes.  A finer grid does not always imply 
better accuracy, but the Grassmann distance criterion correctly identifies the grid size  of $N = 80$ as the one 
generating  the eigenvalue which agrees best best agrees with the value reported in~\cite{Orzag1971}. Note that the decay in accuracy as $N$ grows is due to round-off error from the ill-conditioned Chebyshev differentiation matrices \cite{Bayliss1995}. However, we emphasize that unlike the boundary derivative criterion, the angle criterion is robust in detecting errors due to round-off.

\subsection{Combining Quality Criteria with Implicit Constraints} \label{combination.sec}
\begin{figure}
\centering
\begin{subfigure}[t]{\dimexpr0.30\textwidth+20pt\relax}
    \makebox[20pt]{\raisebox{58pt}{\rotatebox[origin=c]{90}{$N = 16$}}}%
    \includegraphics[width=\dimexpr\linewidth-20pt\relax]
    {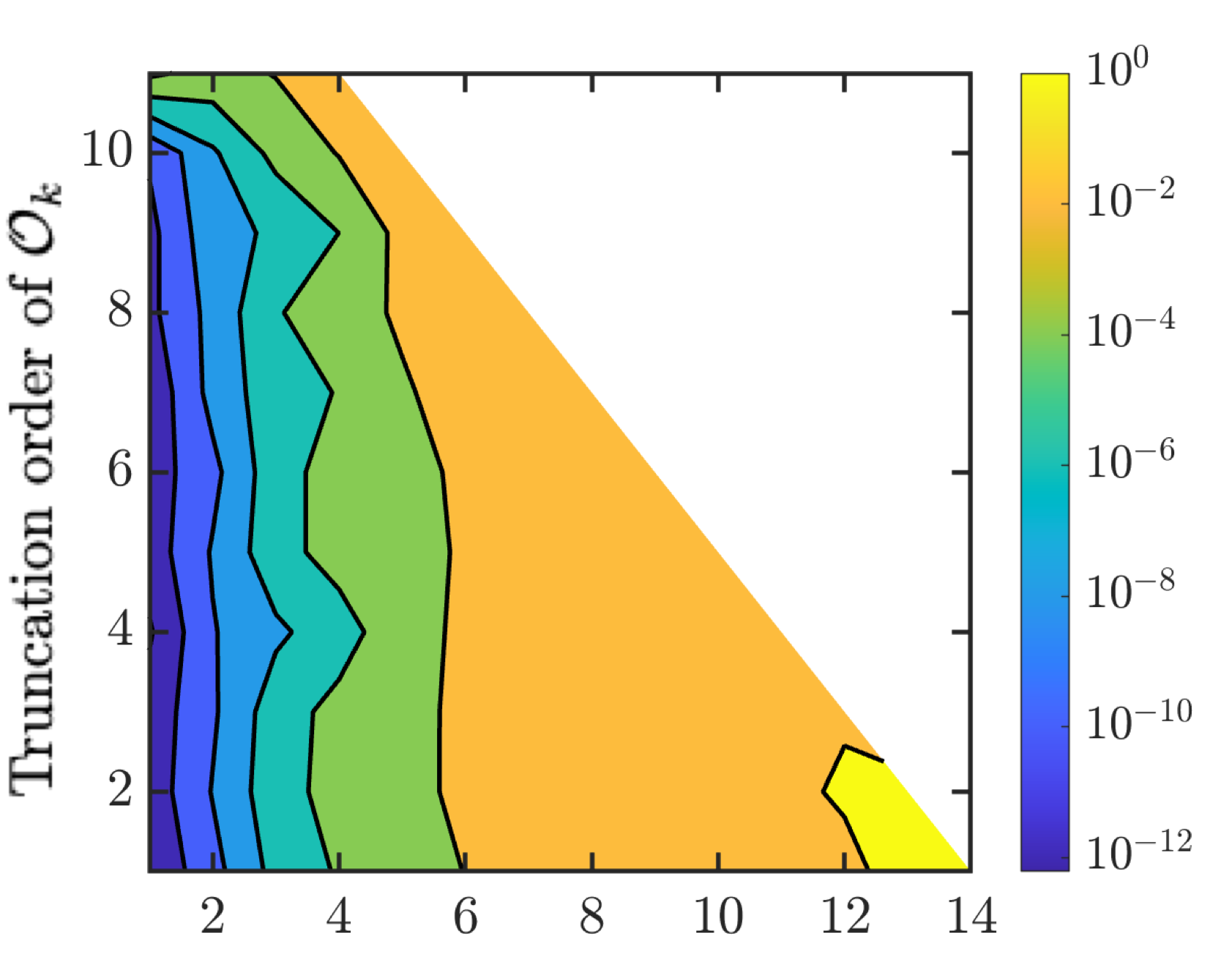}
    \makebox[20pt]{\raisebox{58pt}{\rotatebox[origin=c]{90}{$N = 64$}}}%
    \includegraphics[width=\dimexpr\linewidth-20pt\relax]
    {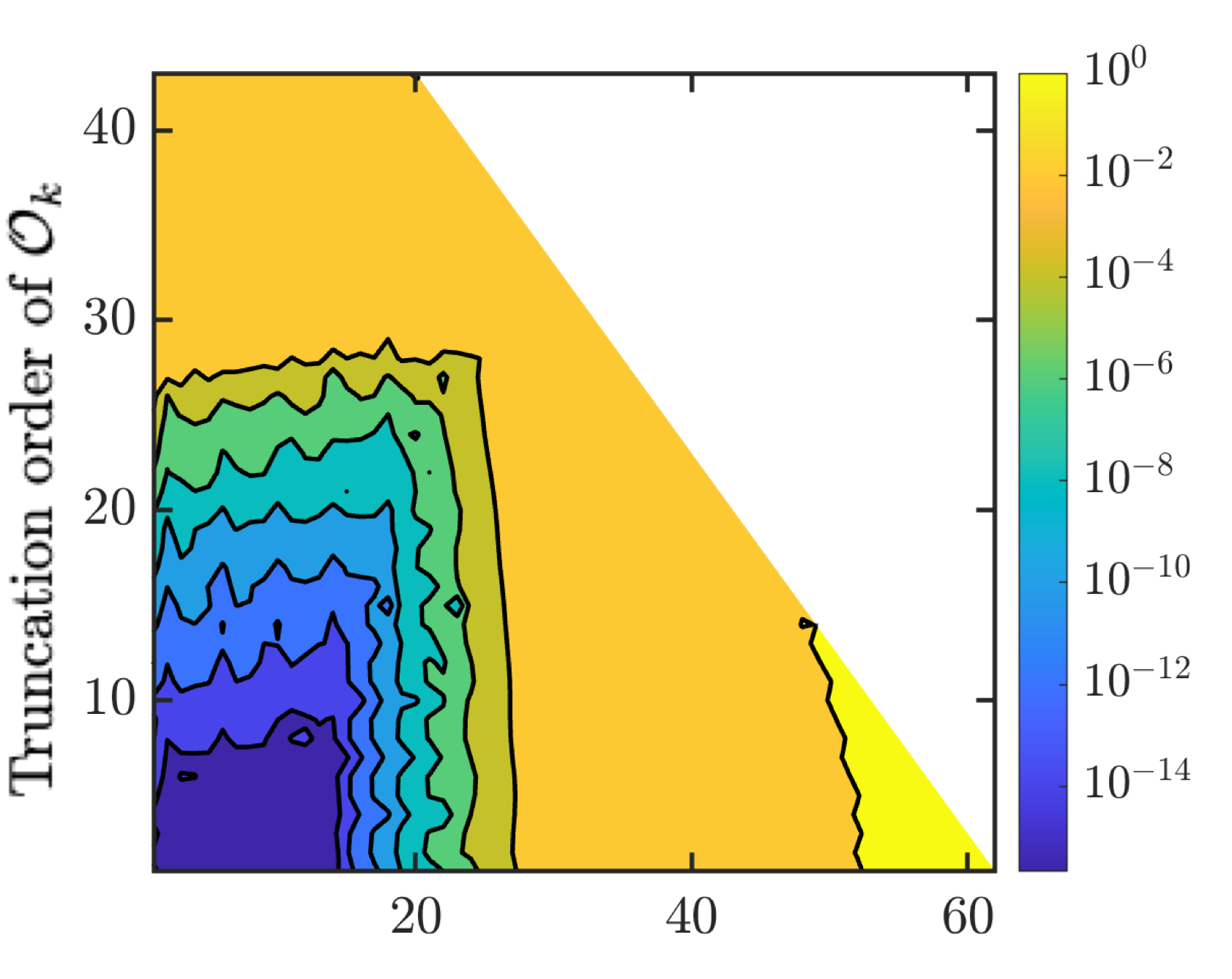}
    \makebox[20pt]{\raisebox{58pt}{\rotatebox[origin=c]{90}{$N = 128$}}}%
    \includegraphics[width=\dimexpr\linewidth-20pt\relax]
    {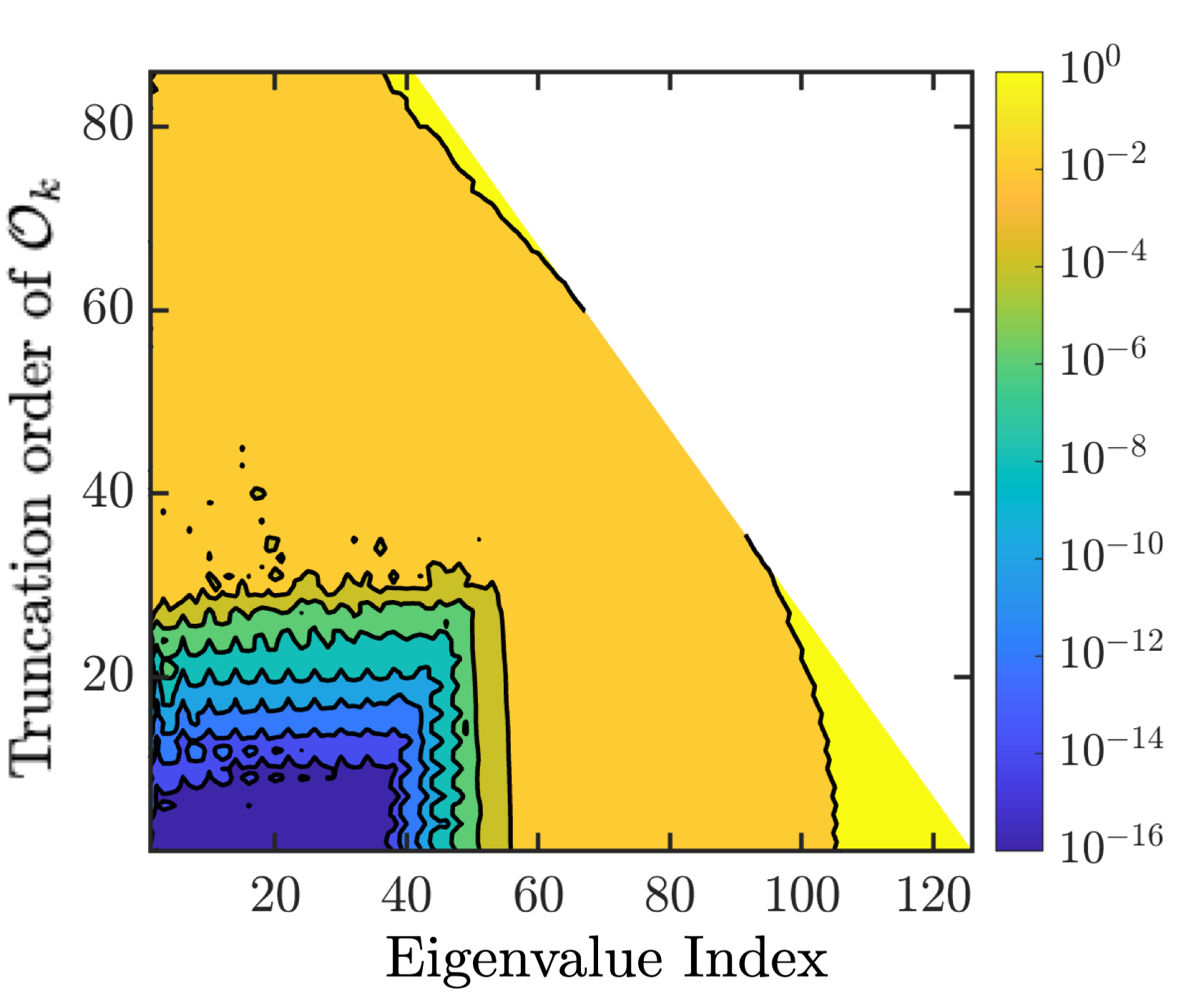}
    \caption{Eigenvalue Error $ | \frac{\bar{\lambda} -\tilde{\lambda}}{\bar{\lambda}}|$ \label{Canuto_1D_contour_eigerror.fig}}
\end{subfigure}\hfill
\begin{subfigure}[t]{0.30\textwidth}
    \includegraphics[width=\textwidth]  
    {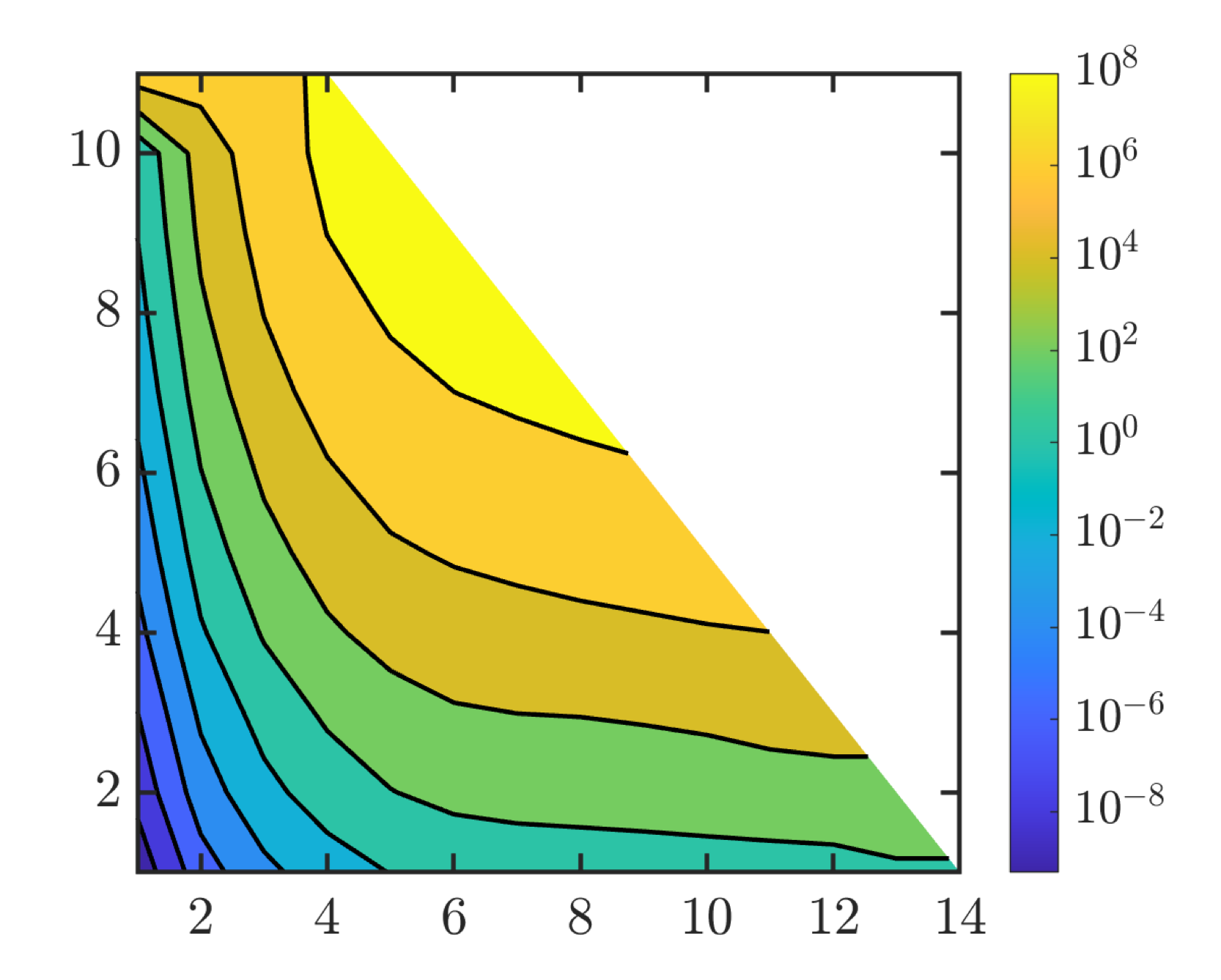}
    \includegraphics[width=\textwidth]
    {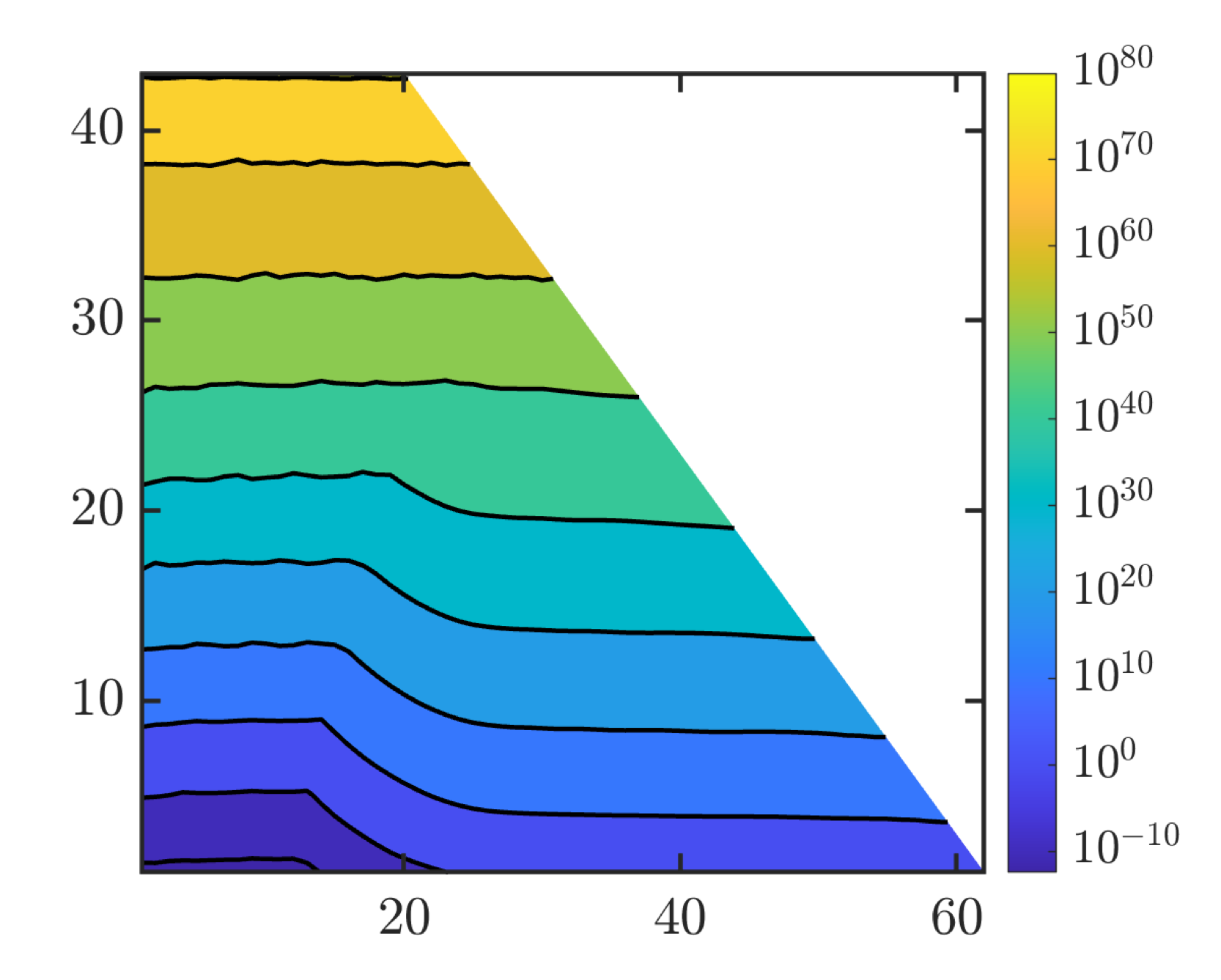}
    \includegraphics[width=\textwidth]
    {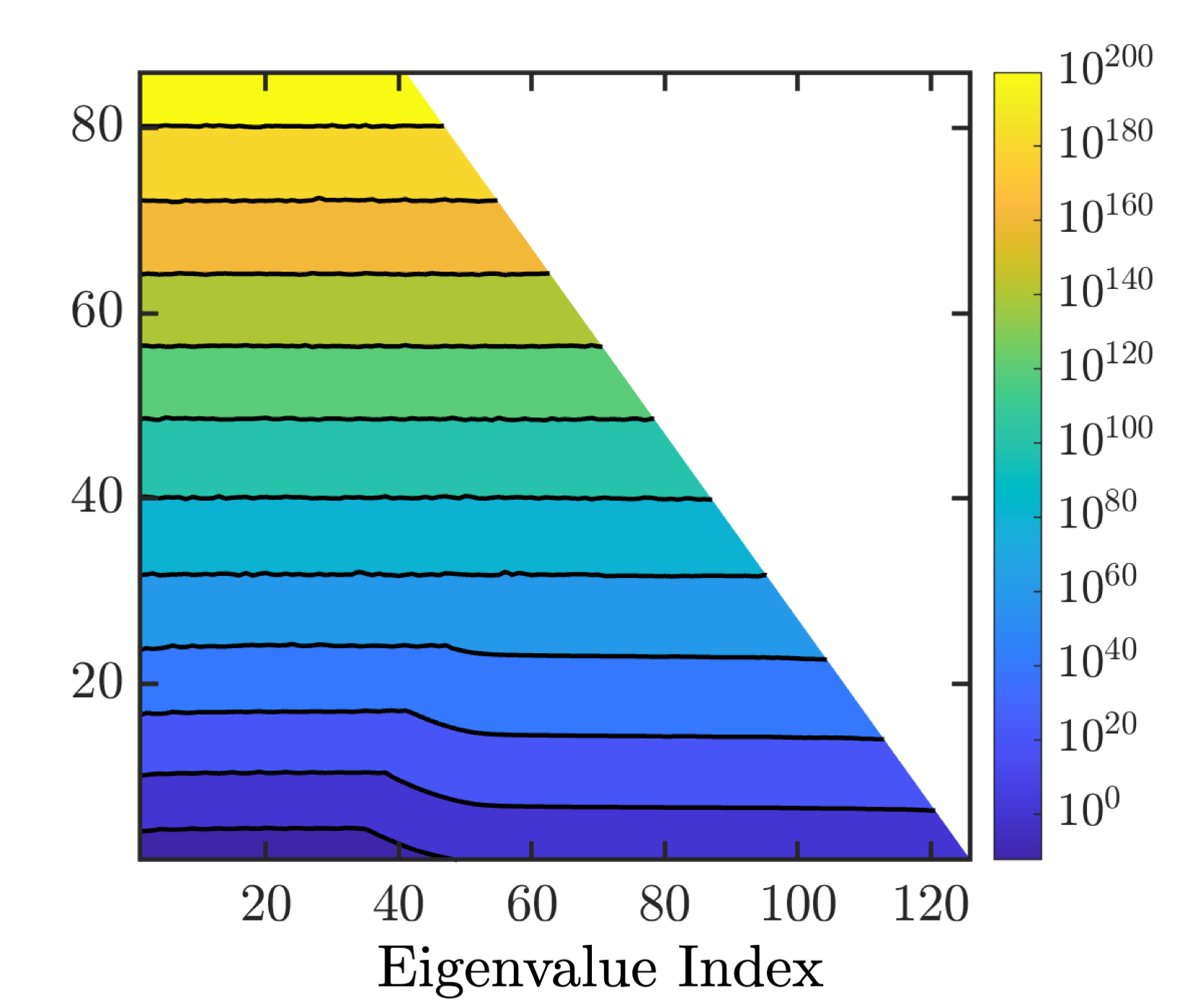}
    \caption{Boundary Derivative $\left\lVert s_{l} \right\rVert$ \label{Canuto_1D_contour_xdot.fig}}
\end{subfigure}\hfill
\begin{subfigure}[t]{0.30\textwidth}
    \includegraphics[width=\textwidth]  
    {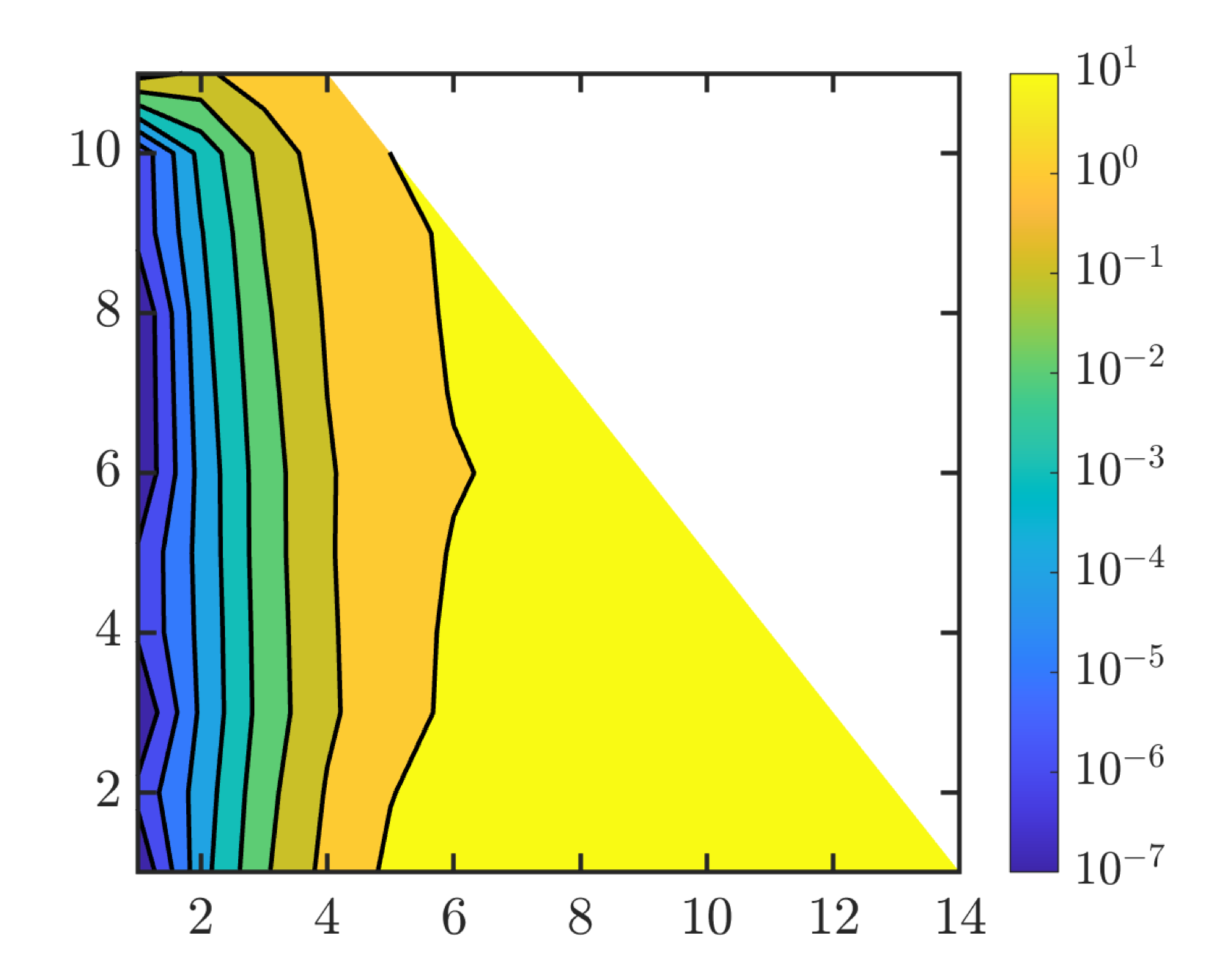}
    \includegraphics[width=\textwidth]
    {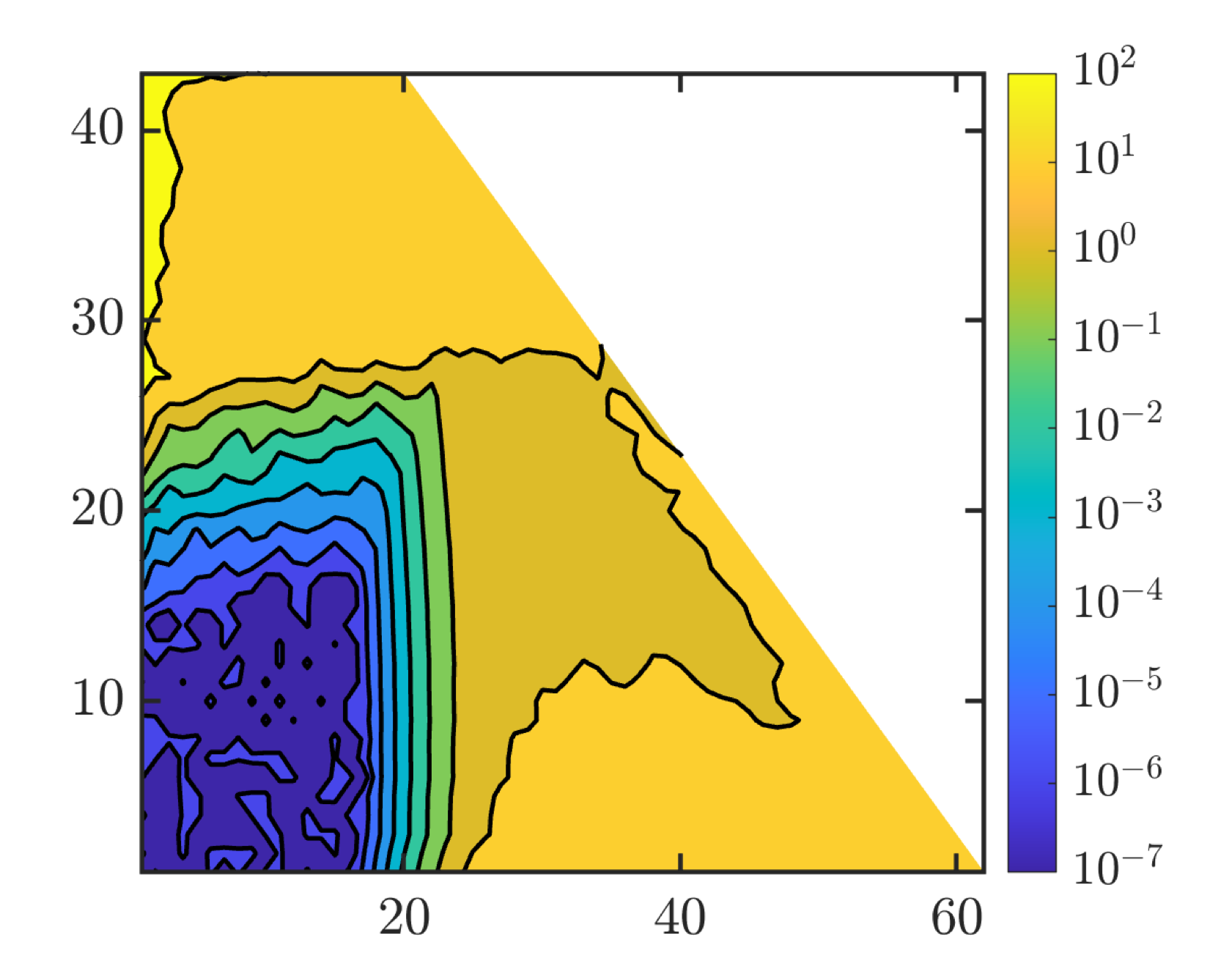}
    \includegraphics[width=\textwidth]
    {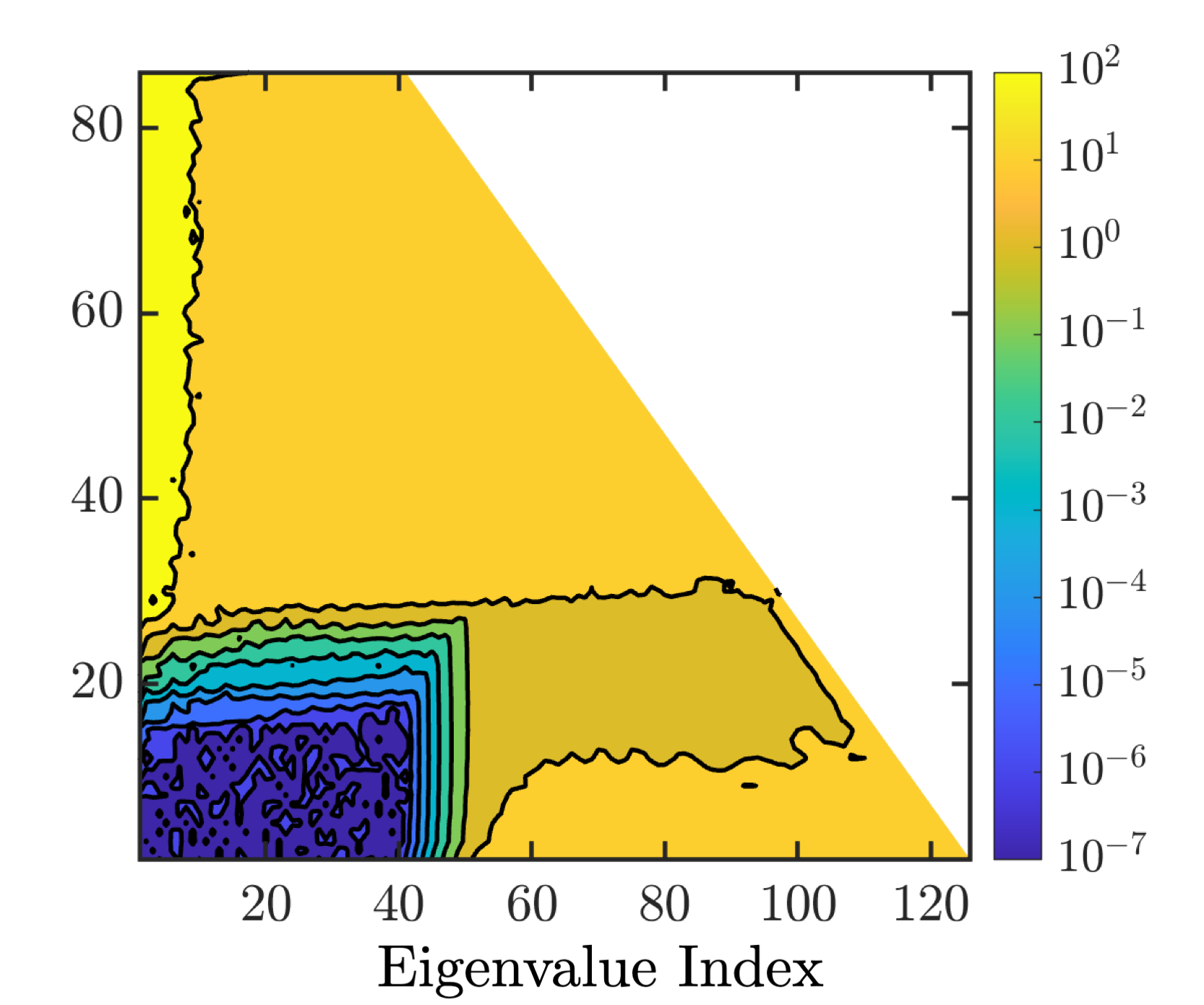}
    \caption{Grassmann Distance $\theta_l$\label{Canuto_1D_contour_grassmann.fig}}
\end{subfigure}
\caption{Contours of the true eigenvalue error (Column (\subref{Canuto_1D_contour_eigerror.fig})), boundary derivative (Column (\subref{Canuto_1D_contour_xdot.fig})) and Grassmann distance (Column (\subref{Canuto_1D_contour_grassmann.fig})) as in system \req{canuto_disc} discretized using a Chebyshev collocation scheme with $N = 16 $ (first row),  $N = 64$ (second row) and  $N = 128$ (third row). For each grid size, level sets of the error criteria are shown for various truncation orders of the $\mathcal{O}^{k}$. This demonstrates that for small $k$, both $ \left\lVert s_l \right\rVert$ and $\theta_l$ both provide good approximations to the eigenvalue error. However, $ \left\lVert s_l \right\rVert$ degrades at higher $k$ (due to the ill-conditioned derivative matrices) and no longer provides insight into the mode quality. Furthermore, increasing $k$ does not improve the subset of the spectrum that is well captured as shown by the nearly vertical contours in both (\subref{Canuto_1D_contour_eigerror.fig}) and (\subref{Canuto_1D_contour_grassmann.fig}) as $k$ is increased. \label{Canuto_contour_all.fig}}
\end{figure}

Before we discuss using the above quality measures as a tool for model reduction, the reader might be curious whether there is any benefit to applying a small number of implicit constraints using Algorithm \ref{Alg_first} combined with the quality measures outlined in the previous section. 
Figure \ref{canuto_all.fig} illustrated that including a small number of implicit constraints generally doesn't degrade the approximation while potentially removing some spurious spectral content. 
However, we find that applying a small number of additional implicit constraints has no impact on the well approximated portion of approximated operator. 
To demonstrate, we briefly return to the discretization of \req{canuto_sys} where we now combine applying additional implicit constraints via algorithm 1 with the quality measures of algorithm 2 and 3.  Figure \ref{Canuto_contour_all.fig} shows the Grassmann distance $\theta_l$ and boundary derivative  $\left\lVert s_l \right\rVert$ plotted alongside the spectrum's relative eigenvalue error. For all three system sizes, we utilize the eigenvalue error in Figure \ref{Canuto_1D_contour_eigerror.fig} as ground truth. Here we find that the application of additional implicit constraints does not improve the dimension of the well approximated subset. After $k > 10$ the quality of the spectrum uniformly begins to degrade. We find that the Grassmann distance qualitatively matches the eigenvalue error well across all three system sizes. Conversely, due to its sensitivity to condition number, the boundary derivative criterion only provides valuable insight for small $k$ and quickly becomes unusable as more constraints are applied. For model reduction purposes we therefore recommend simply applying the boundary conditions using Algorithm \ref{Alg_first} and filtering out the inaccurate portion of the spectrum using the quality measures above. This process is detailed in the section below. 

\section{Model Reduction by Elimination of Low Quality Modes}					\label{model_reduction.sec}

As seen repeatedly in the examples presented in this paper, almost half or two thirds of computed modes are of ``low quality'' as 
quantified by either the boundary derivative or the Grassmann distance criterion. It is natural then to ask how the removal of  those modes effects the accuracy 
of a numerical model. This can be considered as a form of ``model reduction'' if the reduced model has comparable, or possibly better
accuracy than the original high-dimensional model. This kind of model reduction would be particularly useful for example in control systems design. 
The full investigation of these model reduction issues is beyond the scope of the present work, but here we briefly  demonstrate the potential 
with an example.   

%

\begin{figure}[t]
 	\centering
		\begin{subfigure}[t]{.48\textwidth}
			\includegraphics[width=\textwidth]{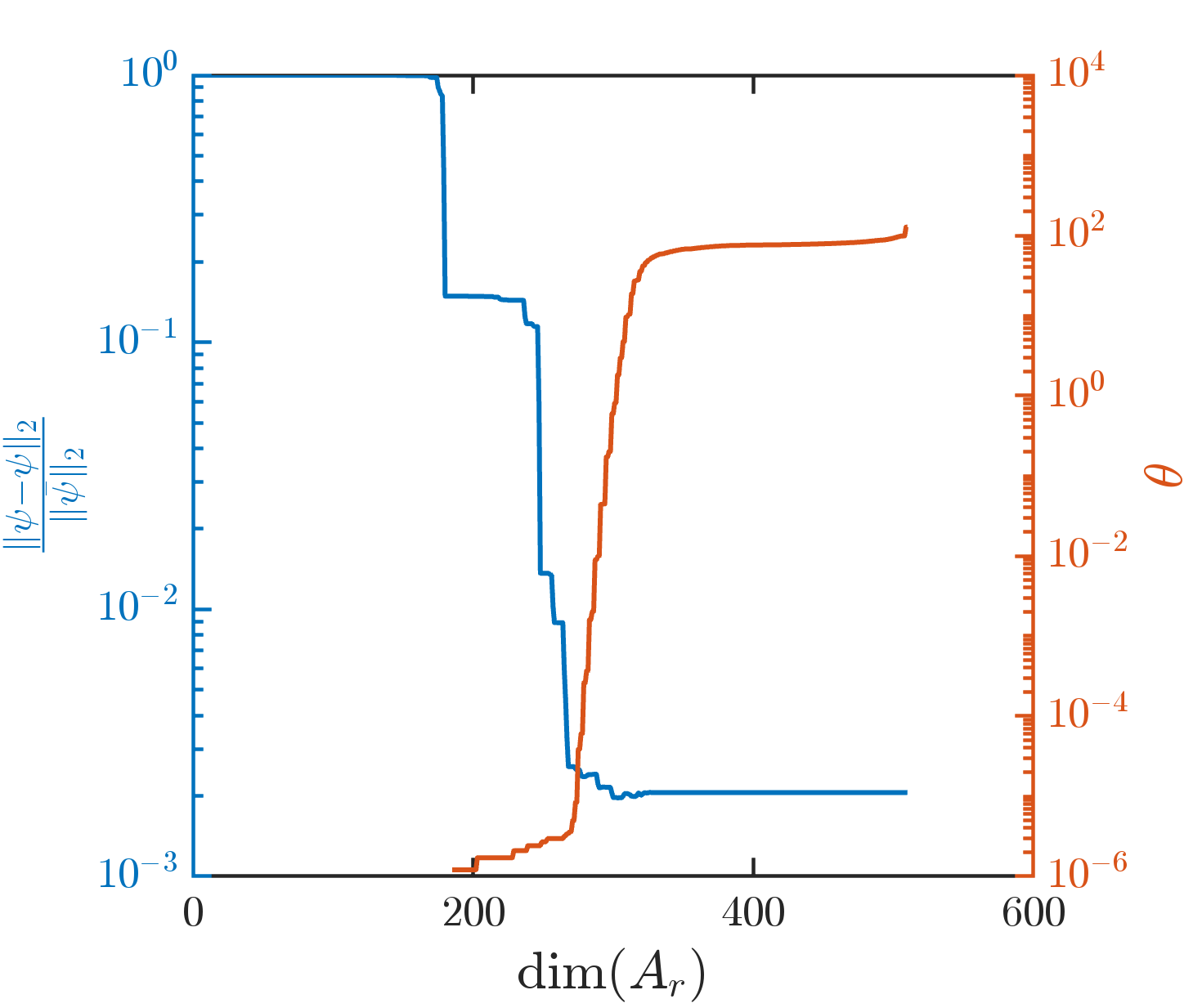}	
			\caption{Bump initial condition \req{acoustic_IC_hard}\label{model_reduc_hard.fig}}
		\end{subfigure}
		\hfill 
		\begin{subfigure}[t]{.48\textwidth}
			\includegraphics[width=\textwidth]{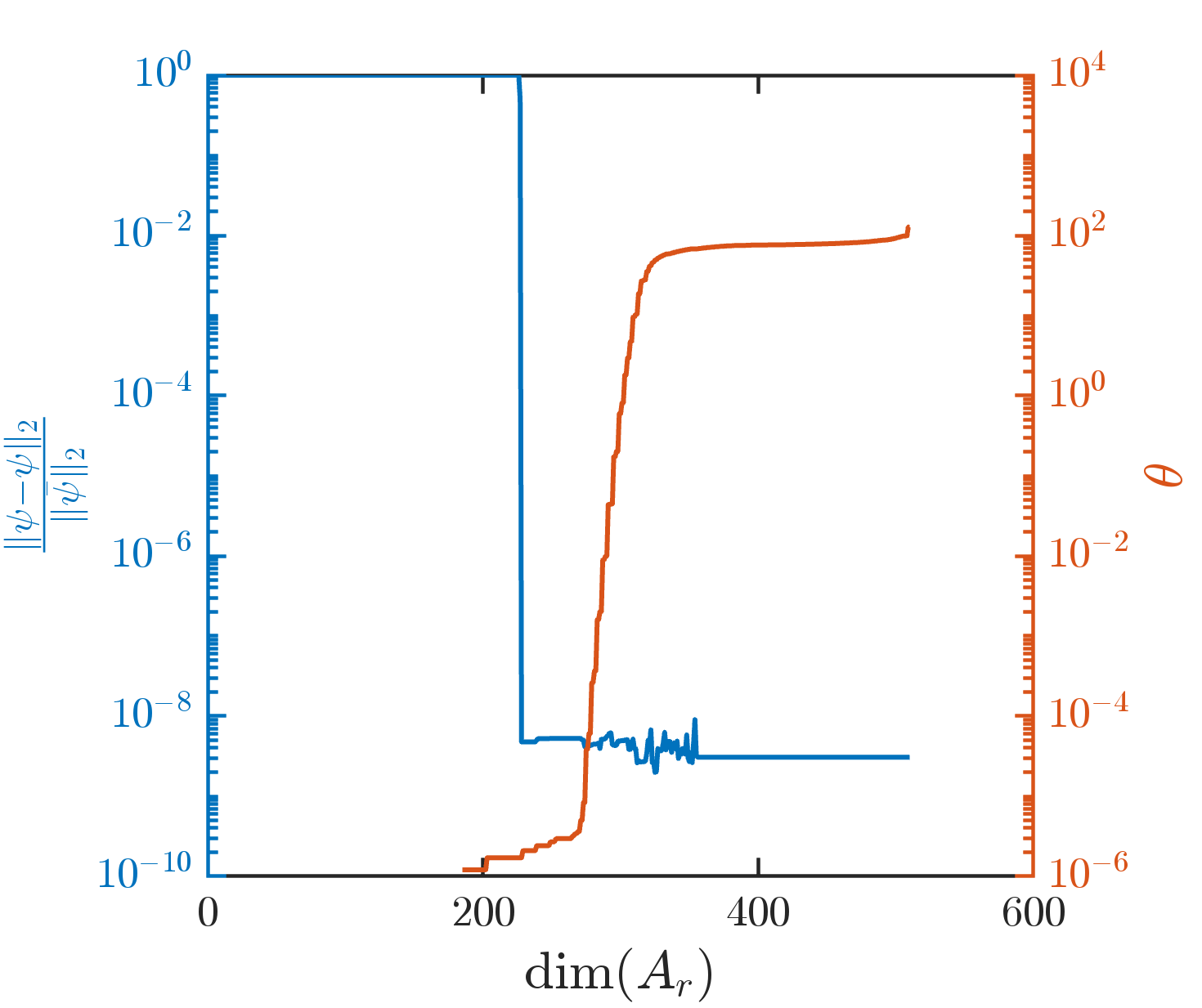}	
			\caption{Sinusoidal initial condition \req{acoustic_IC_easy}\label{model_reduc_easy.fig}}
		\end{subfigure}
		\caption{ Relative error for two different initial conditions of the pressure field of 1-D acoustic wave equation \req{acoustic_sys} at $t = 1$ with $N = 256$ collocation points per field. Relative error is computed using the analytic solution calculated using eigenfunction expansion with 1500 modes. The relative error of the reduced systems $\rm{A_r}$ (left axis) shown against the Grassmann distance of the spectrum of the full-order model sorted by Grassmann distance (right axis). Even in problems such as this where no classically spurious eigenvalues appear, the portion of the spectrum which is identified as poor quality by the Grassmann distance criterion can removed without impacting the quality of the approximation. Each $\rm{A_r}$ is constructed using the $r$ eigenmodes  with smallest Grassmann distance $\theta$, i.e. all eigenmodes depicted to the left of said $\rm{dim(A_r)}$.  Note: for the first $\sim 190$ eigenvalues, $\theta = 0$. \label{model_reduc.fig}}

\end{figure}

Consider the acoustic wave equation in one spatial dimension written in terms of velocity $u$ and pressure $p$ variables with homogeneous 
Dirichlet boundary conditions on pressure  
\begin{align}	
		\partial_t 	\bbm p(x,t) \\  u(x,t)	\\ \ebm  
			&= \bbm 0 & \partial_x \\ \partial_x & 0 \ebm 			\bbm p(x,t) \\ u(x,t) \ebm,
			&  x &\in [-1,1] 			,													\\
		0 &=~ 		p(-1,t) = p(1,t),
		  & t&\geq 0 . \label{acoustic_sys.eq}
\end{align}
The system is discretized with Chebyshev collocation with $N = 256$  collocation points per field, and the boundary conditions are incorporated as in Algorithm~\ref{Alg_first} with $k=1$ resulting in a linear 
vector differential equation 
\be
\dot{x}(t) ~=~ A_{1} ~x(t), 
  \label{full_ord_sys.eq}
\ee
with a $A_1 \in \mathbb{R}^{510 \times 510}$. 
We consider this as the full-order model. To demonstrate the efficacy of using the quality criteria as a model reduction tool, we sort the spectrum of $A_1$ by Grassmann distance (Figure~\ref{model_reduc.fig}).
Then let $\{\rm{A_r}\}$ set of truncated approximates to $A_1$ where each  $A_r$ is constructed using the $r$ modes with smallest Grassmann distance. 
We then simulate each $\rm{A_r}$ forward in time (Matlab: ode45) from the following initial condition conditions
\begin{equation}
	 p(x,0) = \begin{cases}
	 	0, & |x| \ge .3, \\ 
  		e^{-1/\lb 1- \frac{x}{.3})^{4} \rb }, & |x| < .3, 
 \end{cases}
 \hstm \hstm v(x,0) = 0,\label{acoustic_IC_hard.eq} 
\end{equation}
which is particularly challenging for spectral methods since it features sharp fronts as well as a simple sinusoid
\begin{equation}
	p(x,0) = \sin(-\pi x + \pi) 
	\hstm \hstm v(x,0) = 0. \label{acoustic_IC_easy.eq} 
\end{equation}
 As a measure of error, we use the following normalized $L^2[-1,1]$ norm 
\[
	\left.	\left\| \bar{\psi} -  \tilde{\psi} \right\| \right/ \left\| \bar{\psi} \right\|  , 
	\hstm \hstm 
	\psi(x) :=  u(x,1)  ~\mbox{or} ~  p(x,1) , 
\]
where $\bar{\psi}$ is the analytical solution calculated using eigenfunction expansion with 1500 eigemodes, and $\tilde{\psi}$ 
is the solution obtained from the $r$'th order reduced model. The numerical evaluation of the $L^2[-1,1]$ norm
is performed using Clenshaw-Curtis quadrature.  The relative error as a function of $r:={\rm dim}(A_\rmr)$ is shown for each initial condition in Figure~\ref{model_reduc.fig} where it's plotted against the sorted spectrum of \req{full_ord_sys}. Unsurprisingly, the quality of the approximation is dependent on the choice of initial condition. However, due to the nature of the Grassmann distance criteria, it is impossible to discern amongst the well approximated modes (those with Grassmann distance $\mathcal{O}(10^{-6}$ ) which are necessary for a quality approximation. Therefore, even in the sinusoidal initial condition case (Figure \ref{acoustic_IC_easy.eq}) where the first fundamental eigenfrequency is the only necessary mode, the first $\approx 200$ truncated  $\rm{A_r}$ don't contain said mode. Consequently one must always be weary of truncating $\rm{A_r}$ too quickly. 
We therefore recommend visually determining a suitable threshold as the plateauing behavior (Figure \ref{model_reduc.fig} for $\rm{dim}(A_r) \gtrapprox 300$)  is indicative of spurious (or at least inaccurate) spectral content. 
We generally find that choosing a threshold which truncates modes contained in the tail. For the above example a threshold that is $\mathcal{O}(1)$ is suitable. 

\appendix

\section{Proof of Lemma~\ref{DAE2.lemma}} 

	Assume $M$ is a tall,  full column rank $m\times n$ matrix. 
	Let $\ML$ be any left inverse of $M$, and let $N$ be a full column
	 rank $m\times (m\sm n)$ matrix such that $N^*M=0$
	(e.g. columns of $N$ are orthogonal to columns of $M$ and columns of $M$ and $N$ together form 
	a complete basis of $\R^m$) 
	\be
		\tbod{~~N^*~~}{\ML}
		\bbm  \raisebox{.5em}{$M$\rule{0em}{1.5em}} \ebm 
		~=~
		\tbod{0}{I} .
	  \label{NM_I_one.eq}
	\ee	
	Now the matrix on the left is invertible since 
	\begin{align*} 
		0 ~&= 
		\obtd{w_1}{w_2} \tbod{~~N^*~~}{\ML}
		= 
		w_1 N^* + w_2 \ML 
		= 
		w_1 N^*M + w_2 \ML M 
		= 0 + w_2 
		\hsom \Rightarrow\hsom 
		w_2 = 0											\\
		\Rightarrow \hsom 
		0 ~&= w_1 N^* 
		\hstm \Rightarrow \hstm 
		w_1 =0 \hstm \mbox{since $N^*$ has full row rank}. 
	\end{align*} 
	By~\req{NM_I_one}, the inverse has $M$ in its columns and~\req{NM_I_one}  can be completed 
	 by some matrix $Y$ into the inverse relation
	\be
		\tbod{~~N^*~~}{\ML}
		\obtd{ \raisebox{.5em}{$Y$\rule{0em}{1.5em}} }{ \raisebox{.5em}{$M$\rule{0em}{1.5em}} }
		~=~
		\tbtd{I}{0}{0}{I} .
	  \label{NM_I_two.eq}
	\ee

	  We can now rewrite $z(t)$ in terms of 
	two vectors $y(t)$ and $y_\rmc(t)$ as follows 
	\[
		\bbm \raisebox{0.5em}{$z(t)$} \rule{0em}{2em}  \ebm 
		=: 
		\left[ \begin{array}{c:c}  \raisebox{0.5em}{$Y$} \rule{0em}{2em} &  \raisebox{0.5em}{$M$}  \end{array}\right] 
		\left[ \begin{array}{c}  y_\rmc(t) \\ \hdashline y(t) \end{array} \right] , 
		\hspace{4em} 
		\arraycolsep=2pt
		\begin{array}{rcl} 
			\Ims{M} & = &  \Nus{\cO_n} , 					\\
			\R^n &= &  	\Ims{Y} \oplus \Ims{M} . 
		\end{array} 
	\]
	Now  recast the DAE~\req{lin_DAE_2} in terms of the new variables $y$ and $y_\rmc$ 
	\begin{align}
		\arraycolsep=2pt
		\begin{array}{rcl}
			\bbm   \dot{y}_\rmc(t) \\  \dot{y}(t) \ebm 
			 &=& 
			\bbm N^*  \\ \ML	\ebm 			 
			 A 
			\bbm Y  & M	\ebm 
			\bbm   {y}_\rmc(t) \\  {y}(t) \ebm 
				\\ 
            		0 &=& C 
 			\bbm Y  & M	\ebm 
			\bbm   {y}_\rmc(t) \\  {y}(t) \ebm 
		\end{array} 			
		\hstm &\Leftrightarrow \hstm 
		\begin{array}{rcl}
			\bbm   \dot{y}_\rmc(t) \\  \dot{y}(t) \ebm 
			 &=& 
			\bbm N^* A  Y &  N^* A  M 	\\
				 \ML A  Y &  \ML A  M 	\ebm 
			\bbm   {y}_\rmc(t) \\  y(t) \ebm 
				\\ 
            		0 &=& 
 			\bbm C Y  &C M	\ebm 
			\bbm   {y}_\rmc(t) \\  y(t) \ebm 
		\end{array} 										\\		
		\hstm &\Leftrightarrow \hstm 
		\begin{array}{rcl}
			\bbm   \dot{y}_\rmc(t) \\  \dot{y}(t) \ebm 
			 &=& 
			\bbm N^* A  Y &  0 	\\
				 \ML A  Y &  \ML A  M 	\ebm 
			\bbm   {y}_\rmc(t) \\  y(t) \ebm 
				\\ 
            		0 &=& 
 			\bbm C Y  & 0	\ebm 
			\bbm   {y}_\rmc(t) \\  y(t) \ebm 
		\end{array} 						
	  \label{z_trans_struct.eq}
	\end{align}
	The fact that $N^* A  M=0$ follows from recalling that $\Ims{M}$ is $A$-invariant, and therefore 
	$\Ims{AM}\subseteq\Ims{M}$, which is annihilated by $N^*$ by~\req{NM_I_two}. That 
	$CM=0$ follows from recalling from~\req{cO_n_def} that $\Ims{M}=\Nus{\cO_n} \subseteq \Nus{\cO_0}= \Nus{C}$. 

	Finally, observe that the structure of the differential equations implies the following 
	reduction\footnote{In the control theory literature, this reduction is referred to as the ``observability
	decomposition''.}. 
	By Lemma~\ref{DAE1.lemma}, 
	the DAE has a solution iff the initial condition is in $\Nus{\cO_n}$. In this case $y_\rmc(0)=0$, and since 
	$\dot{y}_\rmc(t) = A_\rmc y_\rmc(t)$ (i.e. it is decoupled from $y$), then $y_\rmc(t)=0$ for all $t\geq 0$. 
	Furthermore, because of the zero in the second block of the transformed $C$, there are no 
	constraints on $y$. Therefore, 
	the subsystem $y_\rmc$ and the constraints can be ignored. 
	\hfill $\blacksquare$

\bibliographystyle{IEEEtran} 
\bibliography{JCP_submission} 

\begin{thebibliography}{10}
\providecommand{\url}[1]{#1}
\csname url@samestyle\endcsname
\providecommand{\newblock}{\relax}
\providecommand{\bibinfo}[2]{#2}
\providecommand{\BIBentrySTDinterwordspacing}{\spaceskip=0pt\relax}
\providecommand{\BIBentryALTinterwordstretchfactor}{4}
\providecommand{\BIBentryALTinterwordspacing}{\spaceskip=\fontdimen2\font plus
\BIBentryALTinterwordstretchfactor\fontdimen3\font minus
  \fontdimen4\font\relax}
\providecommand{\BIBforeignlanguage}[2]{{%
\expandafter\ifx\csname l@#1\endcsname\relax
\typeout{** WARNING: IEEEtran.bst: No hyphenation pattern has been}%
\typeout{** loaded for the language `#1'. Using the pattern for}%
\typeout{** the default language instead.}%
\else
\language=\csname l@#1\endcsname
\fi
#2}}
\providecommand{\BIBdecl}{\relax}
\BIBdecl

\bibitem{Bourne2003}
D.~Bourne, ``{Hydrodynamic stability, the Chebyshev tau method and spurious
  eigenvalues},'' \emph{Contin. Mech. Thermodyn.}, vol.~15, no.~6, pp.
  571--579, 2003.

\bibitem{Gardner1989a}
D.~R. Gardner, S.~A. Trogdon, and R.~W. Douglass, ``{A modified tau spectral
  method that eliminates spurious eigenvalues},'' Tech. Rep.~1, 1989.

\bibitem{Huang1994}
W.~Huang and D.~M. Sloan, ``{The pseudospectral method for solving differential
  eigenvalue problems},'' pp. 399--409, 1994.

\bibitem{Lindsay1992}
K.~A. Lindsay and R.~R. Ogden, ``A practical implementation of spectral methods
  resistant to the generation of spurious eigenvalues,'' \emph{International
  Journal for Numerical Methods in Fluids}, vol.~15, pp. 1277--1294, 1992.

\bibitem{Model1990}
G.~B. McFadden, B.~T. Murray, and R.~F. Boisvert, ``{Elimination of spurious
  eigenvalues in the Chebyshev tau spectral method},'' \emph{J. Comput. Phys.},
  vol.~91, no.~1, pp. 228--239, 1990.

\bibitem{Zebib1987}
A.~Zebib, ``{Removal of spurious modes encountered in solving stability
  problems by spectral methods},'' \emph{J. Comput. Phys.}, vol.~70, no.~2, pp.
  521--525, 1987.

\bibitem{Xu1997}
Y.~S. Xu, ``{Origin and elimination of spurious modes in the solution of field
  eigenvalue problems by the method of moments},'' \emph{Asia-Pacific Microw.
  Conf. Proceedings, APMC}, vol.~1, pp. 465--468, 1997.

\bibitem{Straughan1996}
B.~Straughan and D.~W. Walker, ``Two very accurate and efficient methods for
  computing eigenvalues and eigenfunctions in porous convection problems,''
  \emph{Journal of Computational Physics}, vol. 127, pp. 128--141, 1996.

\bibitem{Fornberg1990}
B.~Fornberg, ``An improved pseudospectral method for initial-boundary value
  problems,'' \emph{Journal of Computational Physics}, vol.~91, pp. 381--397,
  1990.

\bibitem{Danilov2019}
\BIBentryALTinterwordspacing
S.~Danilov and A.~Kutsenko, ``{On the geometric origin of spurious waves in
  finite-volume discretizations of shallow water equations on triangular
  meshes},'' \emph{J. Comput. Phys.}, vol. 398, p. 108891, 2019. [Online].
  Available: \url{https://doi.org/10.1016/j.jcp.2019.108891}
\BIBentrySTDinterwordspacing

\bibitem{Schroeder1994}
W.~Schroeder and I.~Wolff, ``{The Origin of Spurious Modes in Numerical
  Solutions of Electromagnetic Field Eigenvalue Problems},'' \emph{IEEE Trans.
  Microw. Theory Tech.}, vol.~42, no.~4, pp. 644--653, 1994.

\bibitem{Manning2007}
\BIBentryALTinterwordspacing
M.~L. Manning, B.~Bamieh, and J.~M. Carlson, ``{Descriptor approach for
  eliminating spurious eigenvalues in hydrodynamic equations},'' vol. 08544,
  pp. 1--13, 2007. [Online]. Available: \url{http://arxiv.org/abs/0705.1542}
\BIBentrySTDinterwordspacing

\bibitem{Phillips1993}
T.~N. Phillips and G.~W. Roberts, ``{The treatment of spurious pressure modes
  in spectral incompressible flow calculations},'' \emph{J. Comput. Phys.},
  vol. 105, no.~1, pp. 150--164, 1993.

\bibitem{Bewley1998}
T.~R. Bewley and S.~Liu, ``Optimal and robust control and estimation of linear
  paths to transition,'' \emph{Journal of Fluid Mechanics}, vol. 365, pp.
  305--349, 1998.

\bibitem{Orzag1971}
S.~A. Orszag, ``{Accurate solution of the Orr–Sommerfeld stability
  equation},'' \emph{J. Fluid Mech.}, vol.~50, no.~4, pp. 689--703, 1971.

\bibitem{canuto1987boundary}
C.~Canuto and A.~Quarteroni, ``On the boundary treatment in spectral methods
  for hyperbolic systems,'' \emph{Journal of Computational Physics}, vol.~71,
  no.~1, pp. 100--110, 1987.

\bibitem{hespanha2018linear}
J.~P. Hespanha, ``Linear systems theory,'' in \emph{Linear Systems
  Theory}.\hskip 1em plus 0.5em minus 0.4em\relax Princeton university press,
  2018.

\bibitem{canuto2007spectral}
C.~Canuto, M.~Y. Hussaini, A.~Quarteroni, and T.~A. Zang, \emph{Spectral
  methods: evolution to complex geometries and applications to fluid
  dynamics}.\hskip 1em plus 0.5em minus 0.4em\relax Springer Science \&
  Business Media, 2007.

\bibitem{Bayliss1995}
A.~Bayliss, A.~Class, and B.~J. Matkowsky, ``Roundoff error in computing
  derivatives using the chebyshev differentiation matrix,'' pp. 380--383, 1995.

\end{thebibliography}

\end{document}